\documentclass[12pt]{amsart}
\usepackage{amssymb}

\newtheorem{theorem}{Theorem}
\newtheorem{example}{Example}
\newtheorem{proposition}{Proposition}
\newtheorem{conjecture}{Conjecture}

\newtheorem{remark}{Remark}
\newtheorem{definition}{Definition}
\newtheorem{lemma}{Lemma}
\newtheorem{problem}{Problem}

\newcommand{\Z}{{\mathbb Z}}

\newcommand{\R}{{\mathbb R}}
\newcommand{\C}{{\mathbb C}}

\newcommand{\RP}{{\mathbb RP}}

\begin{document}

\title[Complements of discriminants of parabolic function singularities]{Complements of discriminants of real parabolic function singularities}

\author{V.A.~Vassiliev}
\address{Weizmann Institute of Science, Rehovot, Israel}
\email{vavassiliev@gmail.com}
\keywords{discriminant, singularity, surgery, morsification}
\subjclass{Primary: 14Q30. Secondary: 14B07, 14P25}

\begin{abstract}
A conjecturally complete list of components of complements of discriminant varieties of parabolic singularities of smooth real functions is given. We also promote a combinatorial program  that enumerates topological types of non-discriminant morsifications of isolated real function singularities and provides a strong invariant of the components of complements of discriminant varieties.
\end{abstract}

\maketitle 

\section{Introduction}

Let $f:(\R^n, 0) \to (\R, 0)$, $d f(0)=0$, be an analytic function singularity, and $F:(\R^n \times \R^l,0) \to (\R,0)$ be an arbitrary its analytic deformation, i.e. a family of functions $f_\lambda \equiv F(\cdot, \lambda): \R^n \to \R,$ $f_0 \equiv f$. The {\it  discriminant} (aka the {\it level bifurcation set}) $\Sigma = \Sigma(F)$ of this deformation is the set of points 
 $\lambda$ from a neighborhood of the origin in the parameter space $\R^l$, such that corresponding functions $f_\lambda$ have real critical points   with zero critical value close to the origin in $\R^n$ ; see e.g. \cite{AVGZ}--\cite{AGLV2}. The discriminant is a subvariety   in $\R^l$ (of codimension 1 in all interesting cases),  which can divide a neighborhood of the origin  into several connected components.  Discriminants appear  in many problems of PDE theory, mathematical physics and integral geometry as singular loci (aka {\it wavefronts}) of important functions defined by integrals depending on parameters, e.g. the fundamental solutions of many  PDEs. The qualitative behavior of these functions in different connected components of complements of discriminant varieties (where they are regular) can be very different; see e.g. \cite{petr}, \cite{Leray}, \cite{Her}, \cite{ABG}, \cite{Gaa}, \cite{AVGZ}, \cite{Var}, \cite{APLT}--\cite{JSing}. Thus, the problem of the enumeration of these components arises as an early step of the study of these functions. 

Discriminants of function singularities also occur in projective geometry as sets projectively dual to neighborhoods of hypersurface flattening points, see \cite{AGLV2}.

The infinite-dimensional space of all small perturbations of an isolated function singularity is well represented by either of certain sufficiently large but finite-parametric families of functions, so-called {\em versal deformations} of this singularity, see e.g. \cite{AVGZ}, \cite{AGLV1}. For simplest singularities of type $A_\mu$, which are represented by functions $x^{\mu+1}$ in one variable, a versal deformation can be chosen in the form \begin{equation}
x^{\mu+1} + \lambda_1 + \lambda_2 x + \lambda_3 x^2 + \dots + \lambda_{\mu} x^{\mu-1},
\label{am}
\end{equation}
where $\lambda_1, \dots , \lambda_\mu$ are the parameters. For the first non-trivial cases $A_2$ and $A_3$, the discriminants in these spaces have respectively the form of a semicubical parabola in $\R^2$ (projectively dual to a neighborhood of an inflection point of a plane curve) and a {\em swallowtail} in $\R^3$, see Fig.~\ref{sc}. In both cases, they divide the spaces of polynomials (\ref{am}) into sets of polynomials with fixed numbers of (simple) roots. 

\unitlength 0.80mm
\linethickness{0.4pt}
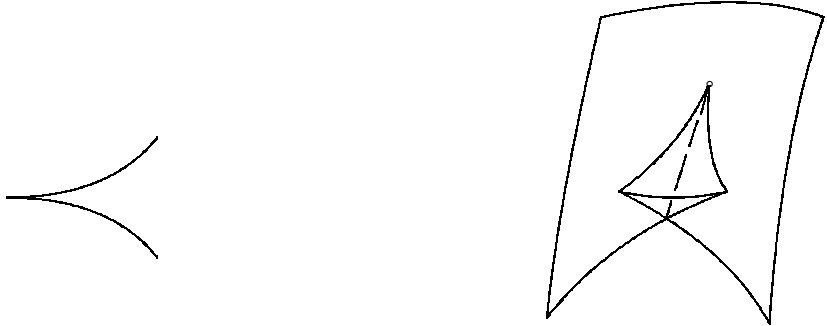
\begin{figure}
\begin{center}
\begin{picture}(30.00,60.00)
\bezier{120}(5.00,24.00)(22.00,24.00)(30.00,34.00)
\bezier{120}(5.00,24.00)(22.00,24.00)(30.00,14.00)
\end{picture}
\begin{picture}(109,60)
\bezier{72}(75.00,25.00)(84.00,23.00)(93.00,25.00)
\bezier{100}(75.00,25.00)(85.00,32.00)(89.80,42.55)
\put(90.10,43.00){\circle{1.00}}
\bezier{76}(90.00,42.55)(89.00,30.00)(93.00,25.00)
\bezier{152}(93.00,25.00)(74.00,18.00)(63.00,4.00)
\bezier{140}(75.00,25.00)(93.00,16.00)(100.00,3.00)
\bezier{208}(100.00,3.00)(102.00,34.00)(109.00,54.00)
\bezier{204}(72.00,54.00)(65.00,24.00)(63.00,4.00)
\bezier{152}(72.00,54.00)(96.00,59.00)(109.00,54.00)
\bezier{12}(82.75,20.45)(83.33,22.00)(83.67,23.33)
\bezier{24}(84.33,25.33)(85.33,28.33)(86.10,31.00)
\bezier{16}(86.67,32.67)(87.00,33.33)(88.00,36.67)
\bezier{15}(89.70,42.00)(89.20,40)(88.60,38)
\end{picture}
\end{center}
\label{sc}
\caption{Discriminants of types $A_2$ (left) and $A_3$ (right)}
\end{figure}

Complements of discriminants and their complex analogs were studied in numerous works,
see e.g. \cite{A70}, \cite{Lo0}, \cite{Loo78}, \cite{Jav}, \cite{Liv}, \cite{Vbook}--\cite{VLoo}, \cite{Ent2}, \cite{sedykh}, \cite{Sed2}. In particular, E.~Looijenga \cite{Loo78} has discovered a natural one-to-one correspondence between the set of connected components of the complement of the discriminant of any {\em simple} (in the sense of V.I.~Arnold) singularity (i.e. a singularity of one of classes $A_k \ (k \geq 1), \pm D_k \ (k \geq 4), \pm E_6, E_7,$ or $E_8$) and certain algebraic structures formulated in the terms of the corresponding Weyl groups. V.D.~Sedykh \cite{sedykh} has enumerated explicitly these components for simple singularities with Milnor numbers at most 6; in \cite{VLoo} this problem was solved for all simple singularities.

The present work is a continuation of  \cite{VLoo}: we now enumerate components of complements of discriminants for the next in difficulty set of singularities, {\em parabolic} ones, see, e.g., \cite{AVGZ}, \cite{LoEll}, \cite{Jav}.  Still, the exposition here is independent on \cite{VLoo}. 
\medskip

There are eight classes of real parabolic singularities: $P_8^1$, $P_8^2$, $+ X_9$, $-X_9$, $X_9^1$, $X_9^2$, $J_{10}^1$, and $J_{10}^3$. 
For our purposes, it suffices to consider singularities of classes $P_8$ of functions of only three variables, and singularities $X_9$ and $J_{10}$ of functions of two variables.
Normal forms of these classes are given in Table \ref{t4}.

\begin{table}
\caption{Real parabolic singularities}
\label{t4}
\begin{center}
\begin{tabular}{|c|l|l|c|c | c|}
\hline
Notation & \qquad  Normal form & Restrictions & $N_{real}$ & $N_{virt}$ & VM \cr
\hline
$P_8^1$ &
$x^3 + \alpha x y z + y^2z + y z^2$ &
$\alpha > -3 $ & 7 & 7 & 6503 \cr $P_8^2$ &
$x^3 + \alpha x y z + y^2z + y z^2$ &
$\alpha < -3 $ & 15 & 15 & 9174 \cr
$\pm X_9$ & $\pm (x^4 + \alpha x^2 y^2 + y^4 )$ & $\alpha^2 < 4$ & 7 & 7 & 27120 \cr
$X_9^1$ & $ ( x^2 - y^2 )(x ^2 + \alpha x y + y^2) $ &
$\alpha^2 < 4 $ & 14 & 10 & 132636 \cr
$X_9^2$ & $ x y (x + y )(y - \alpha x ) $ &
$\alpha > 0$ & 52 & 18 & 134048 \cr
$J_{10}^1$ & $x (x ^2 + \alpha x y^2 + y^4) $ &
$\alpha^2 < 4 $ & 11 & 10 & 2005366 \cr
$J_{10}^3$ & $x (x + y^2)(x - \alpha y^2) $ &
$\alpha >0 $ & 33 & 23 & 2970134 \cr
\hline
\end{tabular} \end{center} \end{table}

\begin{table}
\caption{Versal deformations of parabolic singularities}
\label{vd}
\begin{tabular}{|c|l|}
\hline
Type & Deformation \\
\hline
$P_8$ & $ f_0 + \lambda_1 + \lambda_2 x + \lambda_3 y + \lambda_4 z + \lambda_5 x y + \lambda_6 x z + \lambda_7 y z + \lambda_8 x y z $ \\
$X_9$ & $f_0 + \lambda_1 + \lambda_2 x + \lambda_3 y+ \lambda_4 x^2 + \lambda_5 x y + \lambda_6 y^2 + \lambda_7 x^2 y + \lambda_8 x y^2 + \lambda_9 x^2 y^2 $ \\
$ J_{10} $ & $ f_0 + \lambda_1 + \lambda_2 x + \lambda_3 y+ \lambda_4 x y + \lambda_5 y^2 + \lambda_6 x y^2 + \lambda_7 y^3 + \lambda_8 x y^3 + \lambda_9 y^4 + \lambda_{10} x y^4 $ \\
\hline
\end{tabular}
\end{table}
Standard versal deformations of parabolic singularities are given in Table 2. In any row of this table, $f_0$ is the deformed function given by the corresponding normal form of Table \ref{t4}, and $\lambda_i$ are the parameters of the deformation. 

\begin{theorem}
\label{mtreal}
The numbers of connected components of the complements of the discriminants of versal deformations of parabolic function singularities are not less than the numbers indicated in the fourth column of Table    \ref{t4}.
\end{theorem}
 
\begin{remark} \rm
Conjecturally, these inequalities  are strict, that is, the numbers of components are equal to corresponding numbers from Table \ref{t4}; this conjecture is especially reliable for singularities of types $X_9$ and $J_{10}$.
\end{remark}

We will describe these components in sections \ref{px9}--\ref{P82} devoted to the corresponding singularity types. \medskip

The main tool of our search of components is a combinatorial program enumerating all possible
collections of certain topological invariants of  analytic continuations to $\C^n$ of generic (i.e. strictly Morse and having no zero critical values) perturbations of a singularity.

It is natural to consider  the topological types of such generic complexified perturbations as vertices of a graph, whose edges correspond to generic surgeries of these perturbations, such as Morse surgeries, coincidence of critical values, and jumps of real critical values through zero.
The perturbations $f_\lambda$ having some fixed topological type define a domain in the space $\R^l$ of parameters $\lambda$  of the versal deformation, and their generic surgeries correspond to smooth pieces of the complementary set separating these domains. This graph is connected (as the parameter space is) and finite (since the space of non-generic perturbations of an isolated singularity is diffeomorphic to a semialgebraic set).

Some  edges of this graph are most important for our problem: they correspond to the functions with zero critical value, and hence connect functions from different components of the complement of the discriminant. The enumeration of these components is equivalent to that of the subgraph obtained by removing these edges. 

We replace our graph by its virtual analog, which keeps many its properties and gives us a lot of information on in, including the hints of choosing its vertices (i.e. topological types of real morsifications of our singularity) with prescribed properties.

 Namely, we associate with any non-discriminant strict morsification of a function singularity a collection of discrete topological characteristics, called the {\em virtual morsification}, which is sufficient to predict the result of transformations of this collection under all possible standard topological surgeries of morsifications. 

Given a singularity class, we calculate the virtual morsification related with some its real morsification and run two versions of our program that enumerate  

a) all virtual morsifications that can be obtained from this one by any chains of standard transformations, and 

b) all virtual morsifications that can be obtained from it by chains of transformations modeling the surgeries of functions  that do not change the component of the complement of the discriminant. 

The virtual morsifications found by program b) form a {\em virtual component} of the set of virtual morsifications of our singularity. All real morsifications from a  component of the complement of the (non-virtual) discriminant obviously define virtual morsifications from the same virtual component. Further, we find a virtual morsification found by program a) which does not belong to the list obtained by program b), apply program b)  to it (obtaining a new virtual component), and repeat this procedure until the entire set of virtual morsifications found by program a) is exhausted. 

\begin{theorem}
\label{mtvirt}
The numbers of distinct virtual morsifications and virtual components of parabolic singularities of real functions are as shown in the sixth and fifth columns of Table \ref{t4}, respectively. 
\end{theorem}

After the complete list of virtual components is found, we realize its elements, i.e. find the real morsifications with such topological data, representing corresponding components of the complement of the discriminant variety. \medskip

The relation between the obtained list of virtual components and the set of connected components of the space of non-discriminant perturbations of the initial real function singularity is not very easy, nevertheless, the former list gives a lot of useful information about the latter. Indeed, distinct real components can correspond to the same virtual component, but then these real components are very similar to each other: in all our cases, this situation is related to a certain symmetry of the initial singularity (and therefore of the space of its versal deformation). Also, our program provides rigorous proofs of absence of real morsifications with certain topological properties, and allows us to find the morsifications with some other properties. In addition, it effectively provides a strong invariant of the components of the complement of discriminant: namely, the number of elements of the related virtual component, see \S \ref{invariants}. 
In more details, the algorithm is outlined in \S \ref{outline} 
below and described in chapter V of \cite{APLT}.

\medskip

The case of  parabolic singularities of class $X_9$ is relatively simple: it is closely related to the (well-known) rigid isotopy classification of plane affine curves of degree $4$. The two main differences between these tasks are as follows: 

a) in our case the functions have fixed (or almost fixed) principal homogeneous parts that determine their behavior ``at infinity'', 
and 

b) we consider functions and not their zero sets; 

\noindent
therefore two functions obtained  from one another by a rotation of the plane or multiplication by $-1$ can belong to different components, while the corresponding plane curves are equivalent. 

In a similar way, the enumeration of the components of the space of non-dis\-cri\-mi\-nant perturbations of $P_8$ singularities is related to rigid isotopy classification of smooth and regular at infinity cubic surfaces in $\R^3$. The latter problem is not so simple: an explicit answer to it was published only recently, see \cite{FK}. 
\medskip

\noindent 
{\bf Acknowledgment.} I thank S.~Finashin and V.~Kharlamov very much for writing the article \cite{FK}. A major part of the present research was done during my work at Steklov Mathematical Institute of Russian Academy of Sciences.

\section{Virtual morsifications, their surgeries and algorithms}
\label{outline}

This section contains a reminder of the analogous part of \cite{VLoo} and its adaptation to the case of parabolic singularities.

\subsection{Virtual morsifications and surgeries}
Let $f_0$ be one of polynomials from Table \ref{t4} with  Milnor number $\mu= \mu(f_0)$ (which is equal to the subscript in the notation of its singularity type), 
$\R^\mu$ \ be the parameter space of its versal deformation from Table \ref{vd}, and $f_\lambda$, $\lambda \in \R^\mu$, be a perturbation of $f_0$.

\begin{definition} \rm
Perturbation $f_\lambda$ is {\it generic} if its analytic continuation to $\C^n$ (which we also will denote by $f_\lambda$) has  $\mu$ Morse critical points close to the origin in $\C^n$, and all their critical values are pairwise distinct and not equal to 0.
\end{definition}

Let $f_\lambda$ be generic. Denote by 
$V_\lambda \subset \C^3$ the corresponding {\em Milnor fiber} of function $f_\lambda$ in the case of singularity $P_8$, or of function $f_\lambda+z^2: \C^3 \to \C$ for $X_9$ and $J_{10}$, see \cite{Milnor}, \cite{AVGZ}. (This fiber is defined as  the set $f^{-1}_\lambda(0) \cap D$  (respectively, $(f_\lambda + z^2)^{-1}(0) \cap D$), where $D$ is a small ball centered at the origin in $\C^3$).
 The group $H_{2}( V_\lambda,\Z) $ is then isomorphic to $\Z^\mu$ and is generated by {\em vanishing cycles} defined by a system of {\it non-intersecting} paths connecting the non-critical value $0$ with all these critical values, see, for example, \cite{AVGZ}, \cite{APLT}. The intersection indices of these vanishing cycles in the oriented 4-manifold $V_\lambda$ are well-defined.

\begin{definition} \rm
A system of paths is {\em standard} if 

1) any two paths connecting 0 with complex conjugate critical values are complex conjugate to each other,

2) all paths connecting 0 with real critical values lie (except for their endpoints) in the domain of $\C^1$, where the imaginary part of the coordinate is positive but less than absolute values of the imaginary parts of all non-real critical values of $f_\lambda$.
\end{definition}

\unitlength 1.00mm
\linethickness{0.4pt}
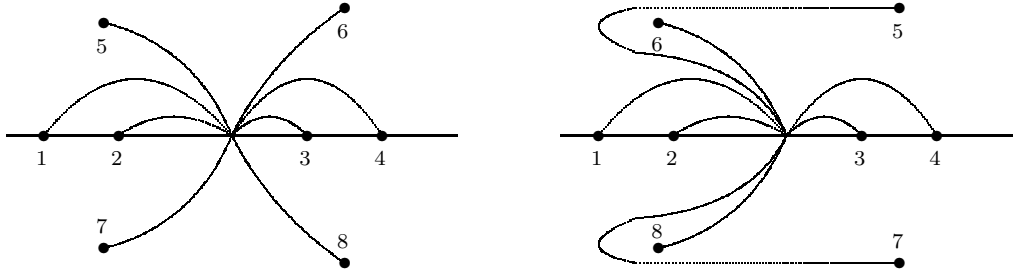
\begin{figure}
\begin{center}
\begin{picture}(60,40)
\put(0,20){\line(1,0){60}}
\put(5,20){\circle*{1.5}}
\put(15,20){\circle*{1.5}}
\put(40,20){\circle*{1.5}}
\put(50,20){\circle*{1.5}}
\bezier{100}(30,20)(17,35)(5,20)
\bezier{50}(30,20)(22,25)(15,20)
\bezier{50}(30,20)(35,25)(40,20)
\bezier{100}(30,20)(40,35)(50,20)
\put(13,35){\circle*{1.5}}
\put(13,5){\circle*{1.5}}
\put(45,37){\circle*{1.5}}
\put(45,3){\circle*{1.5}}
\bezier{100}(30,20)(25,32)(13,35)
\bezier{100}(30,20)(25,8)(13,5)
\bezier{100}(30,20)(35,30)(45,37)
\bezier{100}(30,20)(35,10)(45,3)
\put(4,16){{\tiny 1}}
\put(14,16){{\tiny 2}}
\put(39,16){{\tiny 3}}
\put(49,16){{\tiny 4}}
\put(12,31){{\tiny 5}}
\put(44,33){{\tiny 6}}
\put(12,7){{\tiny 7}}
\put(44,5){{\tiny 8}}
\end{picture} \qquad \quad
\begin{picture}(60,40)
\put(0,20){\line(1,0){60}}
\put(5,20){\circle*{1.5}}
\put(15,20){\circle*{1.5}}
\put(40,20){\circle*{1.5}}
\put(50,20){\circle*{1.5}}
\bezier{100}(30,20)(17,35)(5,20)
\bezier{50}(30,20)(22,25)(15,20)
\bezier{50}(30,20)(35,25)(40,20)
\bezier{100}(30,20)(40,35)(50,20)
\put(13,35){\circle*{1.5}}
\put(13,5){\circle*{1.5}}
\put(45,37){\circle*{1.5}}
\put(45,3){\circle*{1.5}}
\bezier{100}(30,20)(25,32)(13,35)
\bezier{100}(30,20)(25,8)(13,5)
\bezier{100}(30,20)(25,30)(10,31)
\bezier{50}(10,31)(0,35)(10,37)
\bezier{100}(10,37)(35,37)(45,37)
\bezier{100}(30,20)(25,10)(10,9)
\bezier{50}(10,9)(0,5)(10,3)
\bezier{100}(10,3)(35,3)(45,3)
\put(4,16){{\tiny 1}}
\put(14,16){{\tiny 2}}
\put(39,16){{\tiny 3}}
\put(49,16){{\tiny 4}}
\put(12,31.3){{\tiny 6}}
\put(44,33){{\tiny 5}}
\put(12,6.4){{\tiny 8}}
\put(44,5){{\tiny 7}}
\end{picture}
\end{center}
\caption{Standard systems of paths}
\label{standd}
\end{figure}

See Fig.~\ref{standd} for examples of possible standard systems of paths for a real morsification of a singularity with $\mu = 8$.
\medskip

Let us order $\mu$ critical values (and the corresponding vanishing cycles) as follows: first, the real critical values in ascending order, then the values with positive imaginary parts, ordered by the increase of angles at $0$ between the ray of the negative numbers and the corresponding paths, and finally critical values with negative imaginary parts, ordered according to the order (just determined) of their complex conjugate critical values.

The corresponding $\mu$ vanishing cycles can be {\it canonically oriented} by the following conditions:

\begin{enumerate}
\item the canonical orientations of cycles vanishing at real critical points are induced by a fixed orientation of $\R^3$, see \S V.1.6 of \cite{APLT};
\label{ituni}

\item the homology classes of any two vanishing cycles defined by mutually complex conjugate paths are mapped to each other by the complex conjugation in $\C^3$;
\label{itwo}

\item the string of upper triangular matrix elements $a_{1,2}, a_{1,3},$ $ \dots, a_{1,\mu}$, $a_{2,3},$ $ \dots, $ $a_{2,\mu},$ $ \dots, a_{\mu-1,\mu}$ of the matrix of  intersection indices of vanishing cycles is {\em lexicographically maximal} in the class of such matrices defined by systems of orientations satisfying previous two conditions (\ref{ituni}), (\ref{itwo}) (i.e., the number $a_{1,2}$ has the greatest possible value, $a_{1,3}$ has the greatest possible value over the set of matrices with the previous value of $a_{1,2}$, etc.);
\label{itre}

\item if $f_\lambda$ has no real critical points, then the string of $\mu$ intersection indices in $V_\lambda$ of vanishing cycles with the naturally oriented submanifold $V_\lambda \cap \R^3$ is lexicographically maximal in the class of such strings defined by systems of orientations satisfying previous three conditions (\ref{ituni})--(\ref{itre}).
\label{itfou}
\end{enumerate}

Due to the connectedness of Dynkin diagrams of isolated singularities (see, for example, \cite{AVGZ}, vol. 2, \S I.3), conditions (\ref{ituni})--(\ref{itre}) uniquely fix a system of orientations of vanishing cycles in all cases, except for the situation considered in item (4) (when the simultaneous change of all orientations preserves the intersection matrix and complex conjugation). It easily follows from Theorem 1.3 of \cite{APLT}, \S V.1.6, that in the latter case all the intersection indices of vanishing cycles with the set $V_\lambda \cap \R^3$ cannot be simultaneously equal to zero, therefore condition (4) ends with fixing all orientations in this case as well.

\begin{definition} \rm
A {\it virtual morsification} associated with a generic real perturbation $f_\lambda$ is a collection of topological characteristics of $f_\lambda$ (depending also on the choice of a standard system of paths) consisting of 
\begin{itemize}
\item the $\mu \times \mu$ matrix of intersection indices of canonically ordered and oriented vanishing cycles in $V_\lambda$, defined by a standard system of paths, 

\item the string of $\mu$ intersection indices of these vanishing cycles with cycle $V_\lambda \cap \R^3$,

\item the string of Morse indices of all real critical points of $f_\lambda$, ordered by increase of their critical values, and 

\item the number of these points with negative critical values. 
\end{itemize}

Function $f_\lambda$ is called a {\em realization} of any virtual morsification associated with it.

The {\em domain associated with a certain virtual morsification} is the set  of all points $\lambda \in \R^\mu$ having  a neighborhood in $\R^\mu$ such that this virtual morsification is associated with 
all generic functions $f_{\lambda'}$ with $\lambda'$ from this neighborhood. 
\end{definition}

\begin{remark} \rm
1. If $f_\lambda$ has none or one pair of non-real critical points, then there is only one homotopy class of standard systems of paths, and hence only one virtual morsification associated with $f_\lambda$. On the other hand, two systems of paths of the same morsification, shown in the two parts of Fig.~\ref{standd}, define different sets of vanishing cycles and, generally speaking, different matrices of intersection indices.

2. Intersection indices of $V_\lambda \cap \R^3$ with cycles vanishing at real critical points can be deduced from the intersection indices of vanishing cycles with each other and Morse indices or these critical points, see \cite{APLT}, \S V.3.

3 (cf. \cite{Loo78}). The points $\lambda \in \R^\mu$ corresponding to non-generic perturbations $f_\lambda$ having a non-real non-Morse critical point, or a non-real critical point with critical value 0, or a pair of critical points with coinciding non-real critical values, constitute a subset of codimension at least 2 in $\R^\mu$, in particular the sets of such points do not divide any domain into parts. Indeed, a real (i.e. commuting with complex conjugation) function having either of these degeneracies, has also another one, which is complex conjugate to it. Therefore such a function belongs to a stratum of complex codimension at least 2 in the parameter space of the complexification of our versal deformation, and the real codimension of the intersection of such a stratum with the subspace $\R^\mu \subset \C^\mu$ is also at least 2.
\end{remark}

\begin{definition} \rm
{\em Elementary surgeries} of virtual morsifications include three transformations of their data, modeling the generic topological surgeries of the corresponding real morsifications, namely
\begin{itemize}
\item[(s1)] collision of two neighboring real critical values (after which corresponding two critical points either meet and go into the complex domain, or change the order in $\R^1$ of their critical values), 

\item[(s2)] collision of two complex conjugate critical values at a point of the line $\R^1$ (after which corresponding critical points of $f_\lambda$ either meet at a real point and come to the real space, or miss one another, while the imaginary parts of their critical values change their signs), 

\item[(s3)] jumps of real critical values through 0; \\ \hspace*{4cm} and additionally

\item[(s4)] specifically virtual transformations within the classes of virtual morsifications associated with the same real morsification, caused by changes of standard systems of paths going from 0 to imaginary critical values (see Fig.~\ref{standd}).
\end{itemize}

\noindent
We will use notation (s1), (s2) and (s3) for both surgeries of real functions and for transformations of virtual morsifications modeling them.
\end{definition}

\subsection{Case of simple singularities}

 In the case of {\em simple singularities}, any virtual morsification completely determines  the set of surgeries (s1)--(s4) that can be applied to it, as well as the results of all these surgeries; see \S V.8 in \cite{APLT}.

All such possible changes of virtual morsifications of simple singularities can be implemented by real surgeries of functions $f_\lambda$: this follows from the properness of the {\em Looijenga map}. This map acts  from the complex parameter space $\C^\mu$ of the corresponding versal deformation of Table \ref{vd} to the space $ \mbox{Sym}^\mu(\C^1) \simeq \C^\mu$, and sends any parameter value $\lambda$ of the standard complex versal deformation of $f_0$ to the unordered set of critical values of the corresponding function $f_\lambda$ (taken with their multiplicities), see \cite{Lo0}. 

\begin{example} \rm Suppose that a morsification $f_\lambda$
of a simple singularity has critical values shown in Fig.~\ref{standd}.
Let $\zeta \neq 0$ be a real number not separated from 0 by critical values of this  morsification. Consider the path
$\{f_{\lambda(t)}\}$, $t \in [0,1]$, $\lambda(0)=\lambda$, 
 in the corresponding parameter space $\R^8$, starting from the morsification $f_\lambda$ and such that all critical values, except for the two, marked in Fig.~\ref{standd} (left) by numbers 6 and 8, remain fixed, and these two critical values move symmetrically along the segments connecting them with the point $\zeta$. By virtue of the properness of the Looijenga map, such a path exists (and indeed does not leave the space of real perturbations), and its endpoint $f_{\lambda(1)}$ $($i.e. a function with double critical value $\zeta)$ is uniquely determined. The real topological shape of this resulting function is predetermined by the intersection index $\langle \Delta_6, \Delta_8 \rangle$ of vanishing cycles of $f_\lambda$ corresponding to these two critical values. Indeed, since the original singularity was simple, this intersection index can only be 0, 1 or $-1$. If it is equal to 0, then the critical points corresponding to these critical values remain imaginary and separated, and no topological surgery visible in real space occurs. If the intersection index is equal to 1 or $-1$, then two imaginary critical points meet at a real singular point of type $A_2$, which can further split into two real Morse critical points: a minimum and a saddlepoint or a saddlepoint and a maximum depending on the sign of this intersection index.
\end{example}

Our programs ``ignore'' the fact that in the case of non-simple singularities Looijenga map can be non-proper, and produce certain lists of virtual morsifications, whose implementation is not guaranteed a priori (and often is a non-formal problem). Nevertheless, these lists help very much in finding,  enumerating and separating  components of complements of discriminants. Let us describe first the action of our programs in the case of simple singularities.

\subsection{Main and reduced programs}
Let us take a real non-discriminant morsification of a simple singularity, calculate a virtual morsification associated with it (say, using the methods of \cite{gab}, \cite{GZ}, \cite{AC}, \cite{APLT}) and apply to it all possible chains of changes (s1)--(s4). All these transformations over virtual morsifications are defined by easy arithmetic operations, therefore we can realize this process by a computer algorithm, see \cite{pro2}. The result of its work is a complete list of virtual morsifications associated with all non-discriminant real morsifications of the initial singularity. 

If we apply only changes (s1), (s2) and (s4), then  we obtain  a list
containing all virtual morsifications  associated with all real morsifications from a real component of the complement of the discriminant (and, in the case of simple singularities, definitely not containing additional elements). It is natural to call this list a {\em virtual component} of the set of all virtual morsifications of our singularity.

\begin{proposition}
\label{true}
Let $\R^\mu$ be the parameter space of one of the standard versal deformations of simple singularities. If two generic perturbations from different connected components of $\R^\mu \setminus \Sigma$ are associated with the same virtual morsification, then the lists of all virtual morsifications associated with real morsifications from these components completely coincide. 
\end{proposition}

\noindent
Indeed, by the construction any virtual morsification determines its  virtual component. \hfill $\Box$

\begin{definition} \rm
Our combinatorial algorithm, which implements all chains of elementary surgeries and provides all virtual morsifications of a singularity, is called the {\em main program}. The algorithm that implements all chains of elementary changes except for (s3) and provides all virtual morsifications from a single virtual component is called the {\em reduced program}. 
\end{definition}

\begin{remark} \rm In fact, these two programs can be turned into one another by commenting out/uncommenting just two statements.
\end{remark}

\subsection{What fails for not simple singularities and what can be saved}
\label{perils}

The Looijenga map is not proper for non-simple singularities. By this reason, not all chains of virtual surgeries admitted by the previous conditions can actually be realized by paths in the parameter space of a versal deformation.

\begin{remark} \rm
It follows easily from the versality property that Looijenga map of a {\em miniversal} deformation (i.e. a versal deformation depending on the minimal possible number  of parameters, which is equal to $\mu(f_0)$) is a surjection (= local diffeomorphism) on the regular part of parameter space, consisting of strictly Morse perturbations of $f_0$. However, relatively long paths in the configuration space $B(\C^1, \mu)$ generally cannot be lifted to the parameter space of such a deformation (which is a small neighborhood of the origin in $\C^\mu$). 
(According to \cite{Jav}, there exist such liftings to entire space $\C^\mu$ in the case of parabolic singularities and their canonical deformations given in Table \ref{vd}, however, lifting a tiny path, all of whose points are configurations in an arbitrarily small neighborhood of the origin in $\C^1$, we can obtain an arbitrarily long path running along the {\em $\mu$=const} stratum of the singularity, in particular leaving any compact ball in the space of parameters of the versal deformation). 
The more, there is no guarantee that any collision of critical values can be lifted there. 
\end{remark}

Nevertheless we can run our program that creates virtual morsifications without worrying about whether they correspond to real surgeries. For surgeries of types (s1) and (s2) we add the condition that if the intersection index of two vanishing cycles under consideration is not equal to 0 or $\pm 1$, then collision is forbidden (and we make some other natural restrictions). The list of all virtual morsifications obtained by this procedure definitely contains all virtual morsifications associated with  real morsifications of the original singularity. In many cases (including all these considered below) it turns out that all the obtained virtual morsifications are still realizable (although this realization may be an informal additional task). 

If the obtained list of virtual morsifications is finite (which is the case for parabolic singularities), we use the restricted program to enumerate virtual components into which this list splits. Then we try to reconstruct the corresponding real components with these topological data, as well as all the real components obtained from them by the action of the symmetry group of $f_0$, and formulate the conjecture that this is a complete list of real components, cf. Theorem \ref{mtreal}. 

An additional trouble for non-simple singularities is that the analog of Proposition \ref{true} is not proved and probably is wrong for sufficiently degenerate singularities: two virtual morsifications can be connected by an allowed virtual move, the real analog of which cannot be realized to connect the related real morsifications. More precisely, our program assumes that (as in Example 1 above) if the intersection index of vanishing cycles related to neighboring critical values is equal to 0, 1 or $-1$ then a collision of these critical values can be implemented by a movement of functions inside the versal deformation. If this assumption is true for some singularity class, then also the conjecture mentioned in the previous paragraph is true for it. I hope that this assumption holds for parabolic singularities by reasons similar to the argumentation of \cite{Jav}.

\subsection{Does the process converge and what to do otherwise?}
\label{converge}

Problem 1984-9 of \cite{A} (asking whether the number of Dynkin diagrams of a fixed isolated function singularity is finite) was posed by V.~Arnold after discussing the perspectives of (the first version of) this program, first of all the finiteness of lists of virtual morsifications. Although in the case of complex singularities the number of Dynkin diagrams is infinite for all non-simple singularities, in our real problem (when we do only surgeries respecting the complex conjugation) this difficulty is postponed: $X_{10}^1$ is the first singularity type for which the number of obtained virtual morsifications is  probably infinite, see \cite{x10}. Moreover, our program is applicable even in such cases. Indeed,  virtual morsifications with too large coefficients can simply be ignored, since the number of real components is definitely finite, and so any such component can be represented by virtual morsifications with limited coefficients and, moreover, be attained by chains of such limited virtual morsifications. 

\subsection{Two-dimensional features}

In fact, we use two slightly different versions of our programs (both main and reduced): one for singularities of functions of corank 2, that essentially depend on two variables, and the other for all higher dimensions, see \cite{pro2}, \cite{pro}. The reason for this occurs from surgery (s2): in the case of the birth of two real critical points we need to predict their Morse indices. The virtual morsification {\em before} the surgery allows us to predict only the parities of these indices  (which are determined by the sign of the intersection index of the corresponding vanishing cycles). In the case of two variables (when the choice is between  pairs of indices (0,1) and (1,2)) this is enough to determine the indices, but in higher dimensions we generally cannot make correct choice between (0,1) and (2,3) etc. For this reason, in the program handling the singularities of corank $\geq 3$ (including $P_8$), Morse indices of virtual morsifications take values only in $\Z /2\Z$, and not in $\Z$.

\subsection{Invariants of non-discriminant morsifications}
\label{invariants}

Another result of our algorithm is that it defines an effective topological invariant of a 
component of the complement of the discriminant variety: namely, the number of virtual morsifications in the corresponding virtual component.

If this number is infinite, then  the number of virtual morsifications obtained by a program ignoring too complicated virtual morsifications can be used instead, see \S \ref{converge}. More precisely, any convention about how complicated virtual morsifications we agree to ignore defines its own invariant.

 Another invariant of the component of $\R^\mu \setminus \Sigma$   containing a point $\lambda$  is the topological type of the set $f_\lambda^{-1}((-\infty, 0])$. Namely, let $B$ be a small ball centered at the origin in $\R^n$ (where $n=2$ in cases $X_9$ and $J_{10}$, and $n=3$ for $P_8$), and $\partial B$ be its boundary. For all $\lambda$ from a very small neighborhood of the origin of the parameter space $\R^\mu$, the corresponding sets $f_\lambda^{-1}((-\infty, 0]) \cap \partial B$ are topologically the same: these sets form a locally trivial (and hence trivializable) fiber bundle over a neighborhood of the origin in $\R^\mu$. Topological classes of pairs $((f_\lambda^{-1}((-\infty, 0]) \cap B), (f_\lambda^{-1}((-\infty, 0]) \cap \partial B))$ (considered up to homeomorphisms which in the restriction to $f_\lambda^{-1}((-\infty, 0]) \cap \partial B$ are compatible with this trivialization) 
are then constant along the components of the complement of the discriminant, 
and therefore define an invariant of these components. In particular, any homological obstruction to the extension of the homeomorphisms \ $f_\lambda^{-1}((-\infty, 0]) \cap \partial B \to f_{\lambda'}^{-1}((-\infty, 0]) \cap \partial B$ \ for $\lambda \neq \lambda'$, \ defined by this trivialization, to homeomorphisms \ $f_\lambda^{-1}((-\infty, 0]) \cap B \to f_{\lambda'}^{-1}((-\infty, 0]) \cap B$ \ provides an algebraic invariant of the components.

A simplification of this invariant (much weaker, but easy to read from the virtual morsification data) is the Euler characteristic of the relative homology group of this pair, i.e., the difference between the number of real critical values of $f_\lambda$ with a negative critical value and even Morse index and the number of real critical points with a negative value and odd Morse index. 

\begin{definition} \rm 
Two invariants defined above are called {\it cardinality} (notation \ $\mbox{\rm Card}$) and {\it topological type} invariants; the simplification of the latter invariant described in the previous paragraph  is called {\em index} (notation \ $\mbox{\rm Ind}$). 
Numbers $\mbox{\rm Card}$ and $\mbox{\rm Ind}$ will be considered as invariants of virtual components as well.
\end{definition}

\begin{remark} \rm
We will see below that, in contrast to the case of simple singularities, even the virtual component of a generic perturbation $f_\lambda$  generally is not determined by the topological type of the set $f^{-1}_\lambda((-\infty,0])$ (while for simple singularities this topological type determines the real component, see Theorem 2 in \cite{VLoo}).
\end{remark}

\begin{definition} 
\label{scale}
\rm For a given morsification $f_\lambda$ with $\mu$ real critical points and all distinct critical values, its {\em standard scale} is any collection of $\mu+1$ functions of the form $f_\lambda- c_j$, $j=0, 1, \dots, \mu$, where  $\mu+1$ constants $c_j$ represent all intervals into which the critical values of $f_\lambda$ divide the real line.
\end{definition}

Almost all data of virtual morsifications (namely, intersection matrices and Morse indices) 
of functions from the same standard scale are the same, therefore after finding them for some generic morsification  it is convenient to  run our algorithm with these initial data for all elements of  its scale.

\subsection{Topological triviality}

According to \cite{LoEll}, the local topology of the discriminant of a parabolic singularity is constant along any $\mu = \mbox{const}$ stratum. (Actually, the work \cite{LoEll} deals only with complex singularities, but its argumentation implies the same result for the real parts of deformations of singularities  invariant under complex conjugation). Therefore (and since the moduli spaces of all real parabolic singularity classes are path-connected) there is a canonical one-to-one correspondence between the local components of the complements of the discriminants for all parabolic singularities of the same class. In particular, when implementing these components, we can arbitrarily choose the value of the parameter $\alpha$ (see Table \ref{t4}).
\medskip

Each of the next seven sections is devoted to one of the classes of parabolic singularities (not including the class $-X_9$, which can be reduced to $+X_9$ by multiplication by $-1$, preserving all the properties of interest for us now).

\section{$+X_9$: positive definite forms of degree four}
\label{px9}

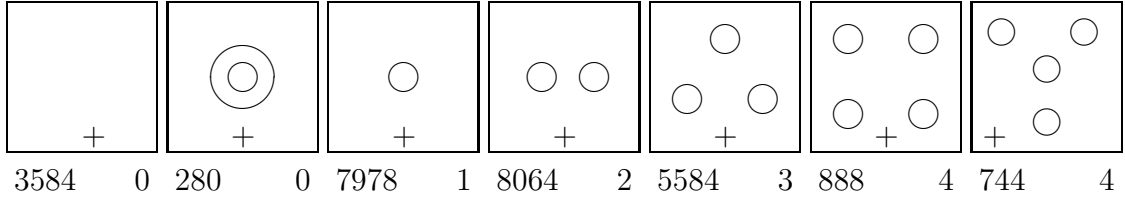
\begin{figure}
\unitlength 1.00mm
\linethickness{0.4pt}
\begin{picture}(20.00,25.00)
\put(0,5){\line(1,0){20}}
\put(0,5){\line(0,1){20}}
\put(0,25){\line(1,0){20}}
\put(20,5){\line(0,1){20}}
\put(10,6){$+$}
 \put(1,0){3584}
\put(17,0){0}
\end{picture}
\begin{picture}(20.00,25.00)
\put(0,5){\line(1,0){20}}
\put(0,5){\line(0,1){20}}
\put(0,25){\line(1,0){20}}
\put(20,5){\line(0,1){20}}
\put(10,15){\circle{4}}
\put(10,15){\circle{8}}
 \put(1,0){280}
\put(17,0){0}
\put(8.5,6){$+$}
\end{picture}
\begin{picture}(20.00,25.00)
\put(0,5){\line(1,0){20}}
\put(0,5){\line(0,1){20}}
\put(0,25){\line(1,0){20}}
\put(20,5){\line(0,1){20}}
\put(10,15){\circle{4}}
 \put(1,0){7976}
\put(17,0){1}
\put(8.5,6){$+$}
\end{picture}
\begin{picture}(20.00,25.00)
\put(0,5){\line(1,0){20}}
\put(0,5){\line(0,1){20}}
\put(0,25){\line(1,0){20}}
\put(20,5){\line(0,1){20}}
\put(7,15){\circle{4}}
\put(14,15){\circle{4}}
\put(8.5,6){$+$}
 \put(1,0){8064}
\put(17,0){2}
\end{picture}
\begin{picture}(20.00,25.00)
\put(0,5){\line(1,0){20}}
\put(0,5){\line(0,1){20}}
\put(0,25){\line(1,0){20}}
\put(20,5){\line(0,1){20}}
\put(5,12){\circle{4}}
\put(15,12){\circle{4}}
\put(10,20){\circle{4}}
 \put(1,0){5584}
\put(17,0){3}
\put(8.5,6){$+$}
\end{picture}
\begin{picture}(20.00,25.00)
\put(0,5){\line(1,0){20}}
\put(0,5){\line(0,1){20}}
\put(0,25){\line(1,0){20}}
\put(20,5){\line(0,1){20}}
\put(5,10){\circle{4}}
\put(5,20){\circle{4}}
\put(15,10){\circle{4}}
\put(15,20){\circle{4}}
\put(8.5,6){$+$}
 \put(1,0){888}
\put(17,0){4}
\end{picture}
\begin{picture}(20.00,25.00)
\put(0,5){\line(1,0){20}}
\put(0,5){\line(0,1){20}}
\put(0,25){\line(1,0){20}}
\put(20,5){\line(0,1){20}}
\put(10,16){\circle{3.5}}
\put(4,21){\circle{3.5}}
\put(10,9){\circle{3.5}}
\put(15,21){\circle{3.5}}
\put(1.5,6){$+$}
 \put(1,0){744}
\put(17,0){4}
\end{picture}
\caption{Perturbations of $+X_9$ singularity}
\label{x9+}
\end{figure}

\begin{proposition}
\label{proX0}
Each singularity of type $+X_9$ has exactly seven virtual components; the topological types of their realizations  are shown in Fig.~\ref{x9+}, and invariants $\mbox{\rm Card}$ and $\mbox{\rm Ind}$ are given respectively by left and right subscripts under these pictures. In particular, the total number of different virtual morsifications of this singularity is equal to 27120.
\end{proposition}

\begin{conjecture}
\label{conX0}
Each of seven virtual components of singularity $+X_9$ is represented by only one component of the complement of the discriminant variety of any its versal deformation.
\end{conjecture}

\noindent
{\it Proof of Proposition \ref{proX0}}. Almost all of these seven pictures (except for the last one) can be realized in the following way. The initial singularity $f_0$ has a perturbation whose zero set is a pair of intersecting ellipses as in Fig. \ref{X901} (left). Let us take a very small further perturbation $\tilde f$ of it which is a strict Morse function. The intersection form of vanishing cycles of the standard system related with $\tilde f$ can be obtained by the method of \cite{GZ}, \cite{AC}. Intersection indices of these cycles with the set of real points can be deduced from them as in \S V.3 of \cite{APLT}, moreover, this calculation is included into our program. Then the ten elements of the standard scale of $\tilde f$ (see Definition \ref{scale}) realize respectively pictures number 1, 3, 4, 5, 6, 5, 4, 3, 2, and 3 of Fig.~\ref{x9+} (counting from the left). Our reduced program shows that corresponding values of invariant \ Card \ are as written under these pictures, and the main program shows that the total number of virtual morsifications is equal to 27120. Also, it is easy to find a virtual morsification not covered by any of these six; substituting it into the reduced program, we get that invariant \ Card \ is equal to 744 for it. Since the sum of all obtained seven values of this invariant is equal to 27120, we conclude that all virtual components are found.

Finally, let us realize the seventh picture of Fig.~\ref{x9+}. The perturbation $ x^4 + \alpha x y + y^4 +\varepsilon y^3$, $ \varepsilon >0$, of singularity $+X_9$ has a singular point of type $+E_6$. Zero level set of this perturbation is an oval with this single singular point. The standard small perturbation of this singularity (described up to a sign in p. 16 of \cite{AC} or in \cite{GZ}, see also Fig.~\ref{e6proof} (left) below) splits this singularity into three Morse saddlepoints with equal critical values and three Morse minima with lower 
values. A level set of obtained function looks as shown in the middle part of Fig.~\ref{X901}. Adding a small positive constant we get the function with zero set of desired shape. 

\unitlength=0.7mm
\begin{figure}
\begin{picture}(40,40)
\put(20,20){\oval(20,40)}
\put(18,3.5){$-$}
\put(18,33.5){$-$}
\put(3,18.5){$-$}
\put(33,18.5){$-$}
\put(18,18.5){$+$}
\put(20,20){\oval(40,20)}
\end{picture} \qquad \qquad \qquad
\unitlength 0.4mm
\begin{picture}(50,70)
\bezier{300}(11,45)(40,75)(47,70)
\bezier{70}(47,70)(50,68)(47,67)
\bezier{250}(47,67)(25,60)(3,67)
\bezier{70}(3,67)(0,68)(3,70)
\bezier{300}(3,70)(10,75)(39,45)
\bezier{500}(39,45)(60,22)(40,6)
\bezier{200}(40,6)(25,-4)(10,6)
\bezier{500}(10,6)(-10,22)(11,45)
\put(22,25){$-$}
\end{picture}
\unitlength 0.7mm
\qquad \qquad \qquad
\begin{picture}(50,40)
\put(28,20){\oval(24,30)}
\put(3,12){\line(1,1){26}}
\put(3,28){\line(1,-1){26}}
\put(11.6,18.8){\small $+$}
\put(26,18.5){\small $-$}
\put(2,18.5){\small $+$}
\put(8.5,30){\small $-$}
\put(8.5,8){\small $-$}
\end{picture}
\caption{Important perturbations for $+X_9$ and $X_9^{1}$}
\label{X901}
\end{figure}
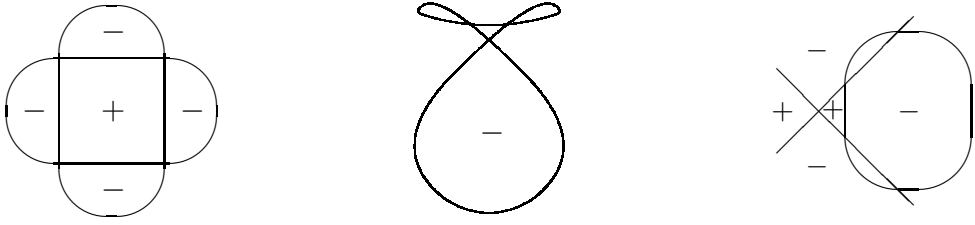

\section{$X_9^1$: two real lines}

There is an involution on the complement of the discriminant variety of $X_9^1$ singularity. If this singularity has the normal form from Table \ref{t4}, and its versal deformation is chosen as in Table \ref{vd}, then this involution sends any function $f_\lambda(x,y)$ to $-f_\lambda(y,x)$. This involution preserves invariant  \ Card. On the other hand,  the sum of values of invariant \ Ind \ of any two components related by this involution is equal to $-1$, in particular this involution acts on the set of components without fixed elements. 

For any value of parameter $\lambda$, the set of negative values $f_\lambda^{-1}((-\infty, 0])$ tends to infinity in $\R^2$ along two opposite sectors, in which $|y|>|x|$. For any non-discriminant perturbation $f_\lambda$ these two sectors can be connected or not connected with each other within the set of negative values of $f_\lambda$. This property holds or fails for all $\lambda$ from any fixed component of $\R^\mu \setminus \Sigma$, so it is a characteristic of this component (and is, of course, a part of the ``topological type'' invariant, see \S \ref{invariants}). Exactly one component from any pair of components related by the above involution  satisfies this property, and the other does not.

\begin{figure}
\unitlength 0.95mm
\linethickness{0.4pt}
\begin{picture}(20.00,25.00)
\put(0,5){\line(1,0){20}}
\put(0,5){\line(0,1){20}}
\put(0,25){\line(1,0){20}}
\put(20,5){\line(0,1){20}}
\put(8,22){\scriptsize $-$}
\bezier{200}(0,6)(10,15)(0,24)
\bezier{200}(20,6)(10,15)(20,24)
 \put(1,0){25874, $-1$}
\end{picture} \ 
\begin{picture}(20.00,25.00)
\put(0,5){\line(1,0){20}}
\put(0,5){\line(0,1){20}}
\put(0,25){\line(1,0){20}}
\put(20,5){\line(0,1){20}}
\bezier{200}(0,6)(10,15)(0,24)
\bezier{200}(20,6)(10,15)(20,24)
\put(10,15){\circle{4}}
\put(8,22){\scriptsize $-$}
 \put(1,0){21200, $-2$}
\end{picture} \ 
\begin{picture}(20.00,25.00)
\put(0,5){\line(1,0){20}}
\put(0,5){\line(0,1){20}}
\put(0,25){\line(1,0){20}}
\put(20,5){\line(0,1){20}}
\bezier{200}(0,6)(10,15)(0,24)
\bezier{200}(20,6)(10,15)(20,24)
\put(10,9){\circle{4}}
\put(10,21){\circle{4}}
\put(3,22){\scriptsize $-$}
 \put(1,0){13212, $-3$}
\end{picture} \
\begin{picture}(20.00,25.00)
\put(0,5){\line(1,0){20}}
\put(0,5){\line(0,1){20}}
\put(0,25){\line(1,0){20}}
\put(20,5){\line(0,1){20}}
\bezier{200}(0,6)(10,15)(0,24)
\bezier{200}(20,6)(10,15)(20,24)
\put(12.5,8){\circle{2.5}}
\put(12.5,22){\circle{2.5}}
\put(9,15){\circle{2.5}}
\put(3,22){\scriptsize $-$}
 \put(2,0){3600, $-4$}
\end{picture} \
\begin{picture}(20.00,25.00)
\put(0,5){\line(1,0){20}}
\put(0,5){\line(0,1){20}}
\put(0,25){\line(1,0){20}}
\put(20,5){\line(0,1){20}}
\bezier{200}(0,6)(10,15)(0,24)
\bezier{200}(20,6)(10,15)(20,24)
\put(8,8){\circle{2.5}}
\put(8,22){\circle{2.5}}
\put(11.5,15){\circle{2.5}}
\put(14,22){\scriptsize $-$}
 \put(2,0){3600, $-4$}
\end{picture} \
\begin{picture}(20.00,25.00)
\put(0,5){\line(1,0){20}}
\put(0,5){\line(0,1){20}}
\put(0,25){\line(1,0){20}}
\put(20,5){\line(0,1){20}}
\bezier{200}(0,6)(13,15)(0,24)
\bezier{200}(20,6)(7,15)(20,24)
\put(17,15){\circle{3}}
\put(8,22){\scriptsize $-$}
 \put(2.5,0){2432, $0$}
\end{picture} \
\begin{picture}(20.00,25.00)
\put(0,5){\line(1,0){20}}
\put(0,5){\line(0,1){20}}
\put(0,25){\line(1,0){20}}
\put(20,5){\line(0,1){20}}
\bezier{200}(0,6)(13,15)(0,24)
\bezier{200}(20,6)(7,15)(20,24)
\put(3,15){\circle{3}}
\put(8,22){\scriptsize $-$}
 \put(2.5,0){2432, $0$}
\end{picture}
\caption{Perturbations of $X_9^1$ singularity}
\label{x91}
\end{figure}
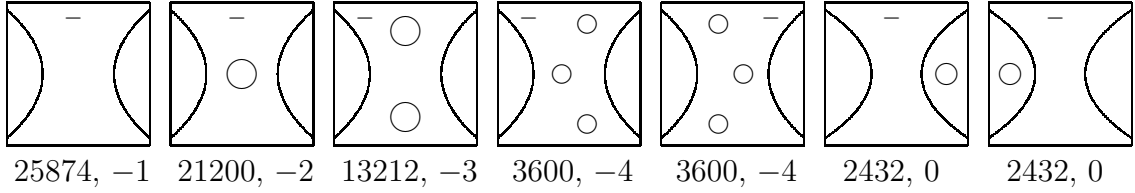

\begin{proposition}
\label{proX1} 1. There are exactly 132636 distinct virtual morsifications of any singularity of type $X_9^1$. 

2. There are exactly ten virtual components of this singularity. They are divided into five pairs; for any pair, the values of invariant \ Card \ are the same, and the sum of values of  invariant \ Ind \ is equal to \ $-1$. 

3. Topological types of realizations of some five virtual components $($repre\-sen\-ting all five pairs$)$ together with the corresponding values of two invariants are given by pictures number 1, 2, 3, 4 and 6  of Fig.~\ref{x91} $($counting from the left$)$. The realizations of the remaining five virtual  components can be obtained from these by the reflection of the pictures in the mirror $\{x=y\}$ and simultaneous change of all signs.

4. Pictures 5 and 7 of Fig.~\ref{x91} represent realizations of the same virtual components as pictures 4 and 6, respectively, but these realizations belong to different components of the complement of the discriminant.
\end{proposition}

\noindent
{\it Proof.}
To realize these pictures by real morsifications, we take first a perturbation of the initial singularity, whose zero set is shown in Fig.~\ref{X901} (right), and then a further minor perturbation $\tilde f$ of this one, which is  a strictly Morse function, such that the critical value at the saddlepoint arising from the crossing point of two lines is greater than the values at four other saddlepoints. Then ten functions $\tilde f - c$ of the standard scale of $\tilde f$ 
 realize pictures 1, 2, 3, 4, 3, 2, 1, 6, and finally the mates of pictures 2 and 1 via our involution.

Virtual morsifications of all elements of this scale and of their mates via this involution
can be calculated by the method of \cite{GZ}, \cite{AC}; substituting them into our programs we obtain that the number of all virtual morsifications of $f_0$ is equal to 132636, and that these elements and their mates represent ten virtual components separated by the pair of invariants $\mbox{\rm Card}$, \ $\mbox{\rm Ind}$. These ten virtual components split into five pairs (corresponding to the orbits of our involution); the values of $\mbox{\rm Card}$ at elements of these pairs are as shown in the left subscripts under the pictures of Fig.~\ref{x91}. The sum of these five values equals exactly a half of 132636, hence our ten virtual components contain all virtual morsifications of $f_0$; this proves statements 1 and 2 of Proposition \ref{proX1}.

Pictures  5 and 7 of Fig.~\ref{x91} can be obtained by a construction symmetric about the vertical axis  to the previous one. Components of \ $\R^9 \setminus \Sigma$ \ described by pictures 6 and 7 have different ``topological type'' invariants; the difference of components described by pictures 4 and 5 follows from Bezout's theorem, since all functions of our versal deformation are polynomials of degree four in $x$ and $y$.
 \hfill $\Box$

\begin{conjecture}
\label{conjX1}
 The complement of the discriminant set in the parameter space of any versal deformation 
of any singularity of type $X_9^1$ contains exactly 14 connected components; their topological types are shown in seven pictures of Fig.~\ref{x91} and in the seven pictures obtained from them by reflection in the diagonal $\{x=y\}$ and change of all signs.
\end{conjecture}

\section{$X_9^2$: four real lines}

In our pictures, we assume that parameter $\alpha$ of the normal form from Table \ref{t4} is equal to $1$, so that the set $f_0^{-1}(0)$ looks like \quad 
\begin{picture}(11,10)
 \put(1.4,1.4){\line(1,1){7.2}}\put(1.4,8.6){\line(1,-1){7.2}}\put(0,5){\line(1,0){10}}\put(5,0){\line(0,1){10}}
\end{picture} .

\unitlength 1.00mm
\linethickness{0.4pt}
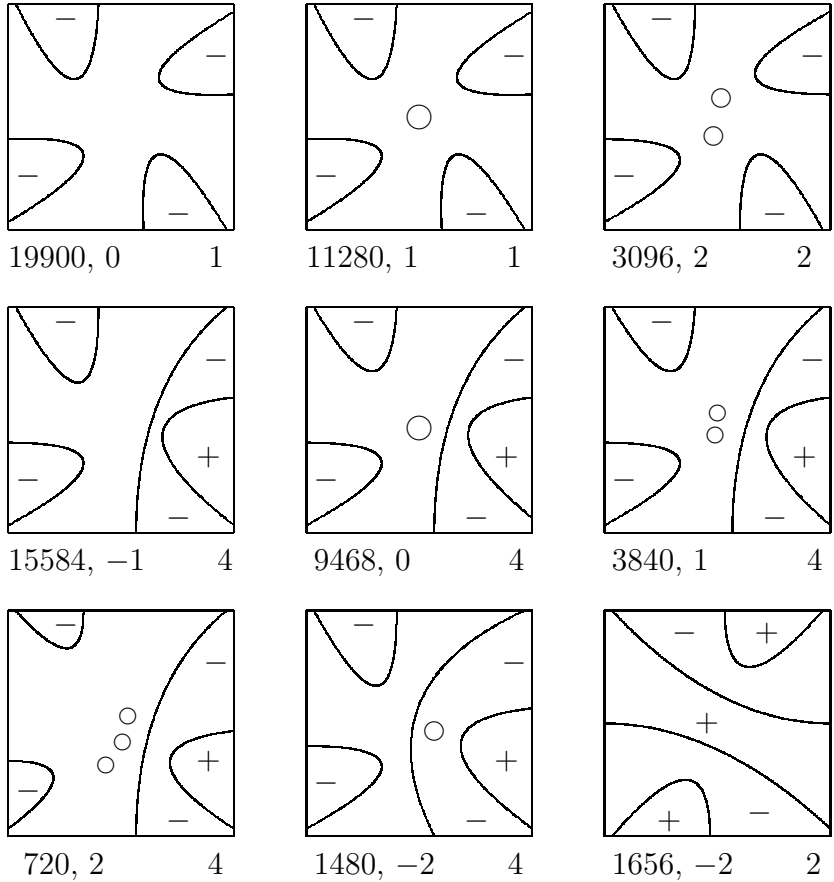
\begin{figure}
\begin{center}
\begin{picture}(30.00,35.00)
\put(0,5){\line(1,0){30}}
\put(0,5){\line(0,1){30}}
\put(0,35){\line(1,0){30}}
\put(30,5){\line(0,1){30}}
\bezier{200}(0,17)(20,18)(0,6)
\bezier{200}(18,5)(17,25)(29,5)
\bezier{200}(30,23)(10,22)(30,34)
\bezier{200}(12,35)(11,15)(1,35)
\put(1,11){$-$}
\put(6,32){$-$}
\put(21,6){$-$}
\put(26,27){$-$}
 \put(0,0){19900, 0 \qquad \, 1}
\end{picture} \qquad 
\begin{picture}(30.00,35.00)
\put(0,5){\line(1,0){30}}
\put(0,5){\line(0,1){30}}
\put(0,35){\line(1,0){30}}
\put(30,5){\line(0,1){30}}
\bezier{200}(0,17)(20,18)(0,6)
\bezier{200}(18,5)(17,25)(29,5)
\bezier{200}(30,23)(10,22)(30,34)
\bezier{200}(12,35)(11,15)(1,35)
\put(15,20){\circle{3}}
\put(1,11){$-$}
\put(6,32){$-$}
\put(21,6){$-$}
\put(26,27){$-$}
 \put(0,0){11280, 1 \qquad \, 1}
\end{picture} \qquad 
\begin{picture}(30.00,35.00)
\put(0,5){\line(1,0){30}}
\put(0,5){\line(0,1){30}}
\put(0,35){\line(1,0){30}}
\put(30,5){\line(0,1){30}}
\bezier{200}(0,17)(20,17)(0,6)
\bezier{200}(18,5)(18,25)(29,5)
\bezier{200}(30,23)(10,23)(30,34)
\bezier{200}(12,35)(12,15)(1,35)
\put(15.5,22.5){\circle{2.5}}
\put(14.5,17.5){\circle{2.5}}
\put(1,11){$-$}
\put(6,32){$-$}
\put(21,6){$-$}
\put(26,27){$-$}
 \put(1,0){3096, 2 \qquad \ \,2}
\end{picture}

\begin{picture}(30.00,40.00)
\put(0,5){\line(1,0){30}}
\put(0,5){\line(0,1){30}}
\put(0,35){\line(1,0){30}}
\put(30,5){\line(0,1){30}}
\bezier{200}(0,17)(20,17)(0,6)
\bezier{200}(12,35)(12,15)(1,35)
\bezier{200}(17,5)(17,25)(29,35)
\bezier{200}(30,23)(11,21)(30,6)
\put(1,11){$-$}
\put(6,32){$-$}
\put(21,6){$-$}
\put(26,27){$-$}
\put(25,14){$+$}
 \put(0,0){15584, $-1$ \qquad 4 }
\end{picture} \qquad 
\begin{picture}(30.00,35.00)
\put(0,5){\line(1,0){30}}
\put(0,5){\line(0,1){30}}
\put(0,35){\line(1,0){30}}
\put(30,5){\line(0,1){30}}
\bezier{200}(0,17)(20,17)(0,6)
\bezier{200}(12,35)(11,18)(1,35)
\bezier{200}(17,5)(17,25)(29,35)
\bezier{200}(30,23)(13,21)(30,7)
\put(15,19){\circle{3}}
\put(1,11){$-$}
\put(6,32){$-$}
\put(21,6){$-$}
\put(26,27){$-$}
\put(25,14){$+$}
 \put(1,0){9468, 0 \qquad \ \, 4}
\end{picture} \qquad 
\begin{picture}(30.00,35.00)
\put(0,5){\line(1,0){30}}
\put(0,5){\line(0,1){30}}
\put(0,35){\line(1,0){30}}
\put(30,5){\line(0,1){30}}
\bezier{200}(0,17)(20,17)(0,6)
\bezier{200}(12,35)(11,18)(1,35)
\bezier{200}(17,5)(17,25)(29,35)
\bezier{200}(30,23)(13,21)(30,6)
\put(15,21){\circle{2}}
\put(14.7,18){\circle{2}}
\put(1,11){$-$}
\put(6,32){$-$}
\put(21,6){$-$}
\put(26,27){$-$}
\put(25,14){$+$}
 \put(1,0){3840, 1 \qquad \ \, 4}
\end{picture}

\begin{picture}(30.00,40.00)
\put(0,5){\line(1,0){30}}
\put(0,5){\line(0,1){30}}
\put(0,35){\line(1,0){30}}
\put(30,5){\line(0,1){30}}
\bezier{200}(0,15)(12,15)(0,6)
\bezier{200}(10,35)(10,25)(1,35)
\bezier{200}(17,5)(17,25)(29,35)
\bezier{200}(30,19)(13,17)(30,6)
\put(15.2,17.5){\circle{2}}
\put(13.0,14.5){\circle{2}}
\put(15.9,21){\circle{2}}
\put(1,10){$-$}
\put(6,32){$-$}
\put(21,6){$-$}
\put(26,27){$-$}
\put(25,14){$+$}
 \put(2,0){720, 2 \qquad \ \ \ 4}
\end{picture} \qquad 
\begin{picture}(30.00,37.00)
\put(0,5){\line(1,0){30}}
\put(0,5){\line(0,1){30}}
\put(0,35){\line(1,0){30}}
\put(30,5){\line(0,1){30}}
\bezier{200}(0,17)(20,17)(0,6)
\bezier{200}(12,35)(12,15)(1,35)
\bezier{200}(17,5)(7,25)(29,35)
\bezier{200}(30,22)(11,20)(30,6)
\put(17,19){\circle{2.5}}
\put(1,11){$-$}
\put(6,32){$-$}
\put(21,6){$-$}
\put(26,27){$-$}
\put(25,14){$+$}
 \put(1,0){1480, $-2$ \qquad 4}
\end{picture} \qquad 
\begin{picture}(30.00,37.00)
\put(0,5){\line(1,0){30}}
\put(0,5){\line(0,1){30}}
\put(0,35){\line(1,0){30}}
\put(30,5){\line(0,1){30}}
\bezier{300}(1,35)(15,20)(30,20)
\bezier{300}(0,20)(15,20)(30,6)
\bezier{200}(14,5)(14,20)(1,5)
\bezier{200}(16,35)(16,20)(29,35)
\put(12,19){$+$}
\put(20,31){$+$}
\put(7,6){$+$}
\put(9,31){$-$}
\put(19,7){$-$}
 \put(1,0){1656, $-2$ \qquad 2}
\end{picture}
\caption{Non-discriminant perturbations of $X_9^2$}
\label{X92}
\end{center}
\end{figure}

We will show that this singularity $f_0$ has, in particular, nine Morse perturbations, topological types of which are drawn in pictures of Fig.~\ref{X92}. In addition, 
the group $\Z_8$ acts on the set of perturbations of $f_0$ and turns the pictures of Fig.~\ref{X92} into other pictures that may belong to different components of the complement of the discriminant. Namely, the {\em negative} of any picture is obtained from it by rotating through the angle $\pi/4$ and changing all signs. Invariant \ Ind \ of any perturbation of $f_0$ is changed by this operation as $I \mapsto -3-I$, so we always get a different component. The square of this operation is just the rotation through the angle $\pi/2$. 

\begin{proposition}
\label{proX2}
1. A singularity of type $X_9^2$ has exactly 134048 distinct virtual morsifications distributed into 18 virtual components. 

2. These 18 virtual components are naturally split into 9 pairs, such that invariant \ Card \ takes equal values on the elements of any pair, and the sum of values of invariant \ Ind \ on these two elements is equal to $-3$.

3. Topological types of some representatives of all nine pairs are shown in Fig.~\ref{X92}, and topological types of their mates in these pairs are represented by the negatives of these pictures; the values of invariants \ Card \ and \ Ind \ of these nine representatives are shown by left-hand subscripts under the corresponding pictures of this figure.

4. For any of nine functions represented in Fig.~\ref{X92}, the number of different components of the complement of the discriminant, visited by the orbit of this function under the action of the group $\Z_4$ generated by the rotation of the plane through $\pi/2$, is at least the right-hand subscript under the corresponding picture; the same numbers hold for the rotations of the negatives of these nine functions.
\end{proposition}

\begin{conjecture} 
\label{conjX2}
The set of components of the complement of the discriminant variety of an arbitrary versal deformation of a singularity of type $X_9^2$ is exhausted by the components represented in Fig.~\ref{X92}, their negatives, and the results of their rotations through the multiples of $\pi/2$: in total 52 elements.
\end{conjecture}

\unitlength=1mm
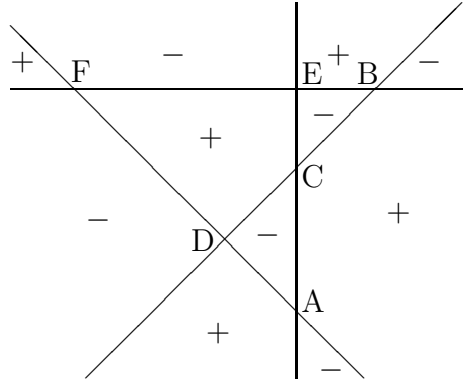
\begin{figure}
\begin{picture}(60,50)
\put(0,38.5){\line(1,0){60}}
\put(38,0){\line(0,1){50}}
\put(10,0){\line(1,1){50}}
\put(47,0){\line(-1,1){47}}
\put(40,34){$-$}
\put(25,31){$+$}
\put(32.5,18){$-$}
\put(26,5){$+$}
\put(50,21){$+$}
\put(10,20){$-$}
\put(54,41){$-$}
\put(0,41){$+$}
\put(20,42){$-$}
\put(41,0){$-$}
\put(42,42){$+$}
\put(38.7,8.8){A} 
\put(46,39.2){B} 
\put(38.7,25.5){C} 
\put(8,39.5){F} 
\put(24,17){D} 
\put(38.7,39.2){E} 
\end{picture} 
\caption{Basic perturbation of $X_9^2$ singularity}
\label{X92sab}
\end{figure}

\noindent
{\it Proof of Proposition \ref{proX2}.} 
First, let us realize all of topological types shown in Fig.~\ref{X92}. 

Take a perturbation of the original singularity of type $X_9^2$, whose set of zeros is shown in Fig.~\ref{X92sab}. Let $\tilde f$ be its tiny perturbation, all critical values of which are different, and the values at the saddlepoints are ordered as $\tilde f(A)<\tilde f(B)<\tilde f(C)<\tilde f(D)<\tilde f(E)<\tilde f(F)$. Then the first six elements of the standard scale of $\tilde f$ realize pictures of Fig.~\ref{X92} with invariant \ Card \ equal to 19900, 11280, 3096, 11280, 19900, 15584.

Second, let $\tilde f_1$ be another minor perturbation of the function from Fig.~\ref{X92sab} such that the critical values at the saddlepoints are ordered as $\tilde f(F)<\tilde f(E)<\tilde f(C)<\tilde f(D)<\tilde f(B)<\tilde f(A)$. The first seven elements of its standard scale provide pictures equivalent up to rotations to pictures with invariant \ Card \ equal to 19900, 11280, 3096, 3840, 9468, 15584, 1480. 

Third, take one more perturbation $\tilde f_2$ of the same function, for which the critical values at saddlepoints $C$ and $F $ are greater than at all other saddlepoints. The function \ $\tilde f_2 - c$ \ from its standard scale, which has exactly six negative critical values, gives the picture with invariant \ Card \ equal to 1656. 

Finally, to get the picture with $\mbox{\rm Card} = 720$ we take the perturbation $f_0 + \varepsilon x^3$ of the initial singularity $f_0$, splitting it to a singularity of type $E_7$ and two separate Morse points, and then perturb it (following \cite{AC}, p. 16) so that the point of type $E_7$ splits into three minima and four saddlepoints, as indicated in the fourth picture of Fig.~6 of \cite{VLoo}, see also Fig.~11 there.

Virtual morsifications of all constructed functions (except for the last one) can be calculated by the method of \cite{GZ}, \cite{AC}. Substituting them into our program, we obtain that the total number of virtual morsifications of our singularity is equal to 134048, and the numbers of virtual morsifications in their virtual components are equal  to leftmost subscripts of corresponding pictures of Fig.~\ref{X92}. Summing up these numbers and taking into account the negatives, we see that there are exactly 1440 virtual morsifications not covered by these 16 virtual components.

Further, the list of all virtual morsifications of our singularity contains one characterized by the intersection matrix
\begin{equation} \left|
\begin{array}{ccccccccc}
$-2$  & $-2$   & 0  &  0 &   1 &   0  &  0  &  1  &  1 \\
   $-2$  & $-2$  &  0  &  0  &  1  &  0  &  0  &  1 &   1 \\
    0  &  0  & $-2$   & 0  &  0  &  1  &  0  &  1  &  0 \\
    0  &  0  &  0  & $-2$  &  1  &  0  &  1  &  0  &  0 \\
    1  &  1 &   0  &  1 &  $-2$  &  0  &  0  &  0  &  0 \\
    0  &  0 &   1  &  0  &  0 &  $-2$  &  0  &  0  &  0 \\
    0   & 0  &  0  &  1  &  0  &  0  & $-2$  &  0 &   0 \\
    1  &  1   & 1  &  0  &  0  &  0 &   0  & $-2$  &  0 \\
    1   & 1  &  0  &  0 &   0  &  0  &  0  &  0  & $-2$ \\
\end{array}
\right| \ , \label{matX9}
\end{equation}
the string of Morse indices of critical points (in ascending order of critical values) equal to (1, 0, 0, 0, 1, 1, 1, 1, 1), and exactly four negative critical values. Our program checks that its virtual component consists of 720 virtual morsifications (and the same is true for its negative). This finishes a proof of statements 1 and 2 of our proposition. 

Statement 4 for all pictures of Fig.~\ref{X92}, except for the third one (with invariant \ Card \ equal to 3096),  follows from the ``topological type'' invariant from \S \ref{invariants}. For the remaining picture, it follows from Bezout's theorem. \hfill $\Box$

\section{$J_{10}^1$ }

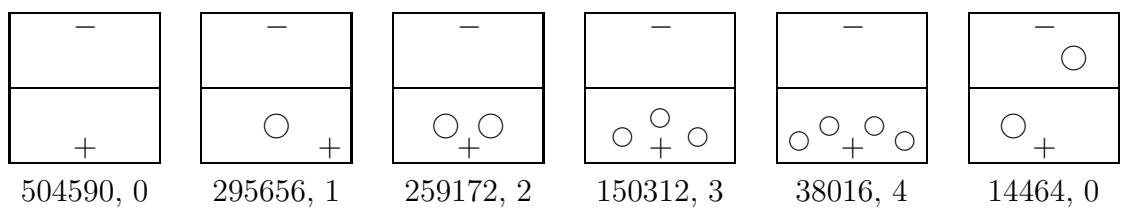
\begin{figure}
\begin{picture}(20.00,25.00)
\put(0,5){\line(1,0){20}}
\put(0,5){\line(0,1){20}}
\put(0,25){\line(1,0){20}}
\put(20,5){\line(0,1){20}}
\put(0,15){\line(1,0){20}}
\put(8.5,22){$-$}
\put(8.5,6){$+$}
 \put(1.5,0){504590, 0}
\end{picture} \quad
\begin{picture}(20.00,25.00)
\put(0,5){\line(1,0){20}}
\put(0,5){\line(0,1){20}}
\put(0,25){\line(1,0){20}}
\put(20,5){\line(0,1){20}}
\put(0,15){\line(1,0){20}}
\put(8.5,22){$-$}
\put(15.5,6){$+$}
\put(10,10){\circle{3}}
 \put(1.5,0){295656, 1}
\end{picture} \quad
\begin{picture}(20.00,25.00)
\put(0,5){\line(1,0){20}}
\put(0,5){\line(0,1){20}}
\put(0,25){\line(1,0){20}}
\put(20,5){\line(0,1){20}}
\put(0,15){\line(1,0){20}}
\put(7,10){\circle{3}}
\put(13,10){\circle{3}}
\put(8.5,22){$-$}
\put(8.5,6){$+$}
 \put(1.5,0){259172, 2}
\end{picture} \quad
\begin{picture}(20.00,25.00)
\put(0,5){\line(1,0){20}}
\put(0,5){\line(0,1){20}}
\put(0,25){\line(1,0){20}}
\put(20,5){\line(0,1){20}}
\put(0,15){\line(1,0){20}}
\put(5,8.5){\circle{2.5}}
\put(15,8.5){\circle{2.5}}
\put(10,11){\circle{2.5}}
\put(8.5,22){$-$}
\put(8.5,6){$+$}
 \put(1.5,0){150312, 3}
\end{picture} \quad
\begin{picture}(20.00,25.00)
\put(0,5){\line(1,0){20}}
\put(0,5){\line(0,1){20}}
\put(0,25){\line(1,0){20}}
\put(20,5){\line(0,1){20}}
\put(0,15){\line(1,0){20}}
\put(3,8){\circle{2.5}}
\put(17,8){\circle{2.5}}
\put(7,10){\circle{2.5}}
\put(13,10){\circle{2.5}}
\put(8.5,22){$-$}
\put(8.5,6){$+$}
 \put(2.5,0){38016, 4}
\end{picture} \quad
\begin{picture}(20.00,27.00)
\put(0,5){\line(1,0){20}}
\put(0,5){\line(0,1){20}}
\put(0,25){\line(1,0){20}}
\put(20,5){\line(0,1){20}}
\put(0,15){\line(1,0){20}}
\put(6,10){\circle{3}}
\put(14,19){\circle{3}}
\put(8.5,22){$-$}
\put(8.5,6){$+$}
 \put(2.5,0){14464, 0}
\end{picture}
\caption{Topological types of perturbations of $J_{10}^1$}
\label{J101}
\end{figure}

\begin{proposition}
\label{proJ1}
1. Any singularity of type $J_{10}^1$ has exactly 2005366 different virtual morsifications divided into ten virtual components.

2. There is an involution on the set of these virtual components, which preserves invariant \ Card \ and replaces the value of invariant \ Ind \ by its opposite. Two virtual components with \ $\mbox{\rm Ind} = 0$ are invariant under this involution, and eight others $($with absolute values of invariant $\mbox{\rm Ind}$ equal to 1, 2, 3 and 4$)$ are split into pairs.

3. Topological types of some realizations of some representatives of all six orbits of this involution are shown in six pictures of Fig.~\ref{J101}, and the values of their two invariants are indicated under these pictures. Topological types of realizations of some representatives of remaining four virtual components can be obtained from four pictures of this figure with \ $\mbox{\rm Ind} \neq 0$ \ by the reflection in a horizontal mirror and simultaneous change of all signs.

4. The rightmost picture of Fig.~\ref{J101} and its reflection in the horizontal mirror describe two realizations of the same virtual morsification, however these realizations belong to different components of the complement of the discriminant.
\end{proposition}

\begin{conjecture}
\label{conjJ1}
There are exactly 11 components of the complement of the discriminant variety of any singularity of type $J_{10}^1$; the topological types of their representatives are as in six pictures of Fig.~\ref{J101} and five pictures, obtained from these pictures, except for the leftmost one, by the reflection in a horizontal mirror and simultaneous change of all signs.
\end{conjecture}

\unitlength=0.8mm
\begin{figure}
\begin{picture}(120,55)
\bezier{300}(60,30)(30,40)(20,40)
\bezier{300}(60,30)(30,20)(20,20)
\bezier{200}(20,40)(0,40)(0,30)
\bezier{200}(20,20)(0,20)(0,30)
\bezier{300}(60,30)(90,40)(100,40)
\bezier{300}(60,30)(90,20)(100,20)
\bezier{200}(100,40)(120,40)(120,30)
\bezier{200}(100,20)(120,20)(120,30)
\bezier{300}(60,55)(80,0)(90,0)
\bezier{300}(120,55)(100,0)(90,0)
\put(20,29){$-$}
\put(65,29){$-$}
\put(114,29){$-$}
\put(90,29){$+$}
\put(90,12){$-$}
\end{picture}
\caption{Starting perturbation for singularity $J_{10}^1$}
\label{j101}
\end{figure}
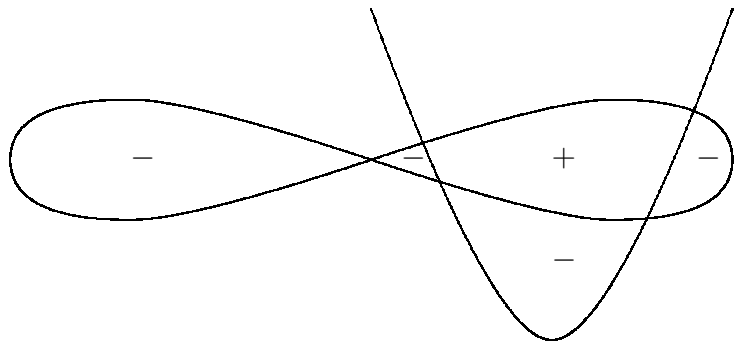

\noindent
{\it Proof of Proposition \ref{proJ1}.}
It is easy to construct a perturbation of a singularity of type $J^1_{10}$, whose zero level set is topologically situated as shown in Fig.~\ref{j101} (where the axis $\{x=0\}$ is assumed to be  horizontal). Let $\tilde f$ be a small further perturbation, all of whose critical values are different, and the critical value at the central saddlepoint (i.e., at the self-intersection point of the \ $\infty$-like curve) is greater than the values at all other saddlepoints. The first (second, third, fourth, fifth, ninth, respectively) element \ $\tilde f - c$ \ of the standard scale of \ $\tilde f$ \ realizes then picture No. 1 (respectively, 2, 3, 4, 5, 6) of Fig. \ref{J101}. 

The entire this construction can be reflected in the horizontal axis; this gives us five new pictures.

The number of all virtual morsifications of our singularity and the numbers of them in all ten virtual components (i.e., the values of invariant  \ $\mbox{\rm Card}$)  are found by our program. The sum of the latter ten numbers is equal to the total number of virtual morsifications of our singularity, hence there are no other virtual components of singularity $J_{10}^1$.

The difference between the  real components represented by the rightmost picture of Fig.~\ref{J101} and its mirror image follows from Bezout's theorem: the restriction of  $f_\lambda$ to any vertical line is a polynomial of degree 3, so two ovals cannot pass one above the other.

\section{$J_{10}^3$} 

For simplicity we will consider the function \ $f_0$ \ with singularity of this type, having the  module \ $\alpha$ \ of the normal form (see Table \ref{t4}) equal to 1; so, zero level set of \ $f_0$ \ looks like \
\begin{picture}(18,8)
\put(0,4){\line(1,0){16}}
\bezier{100}(0,0)(8,8)(16,0)
\bezier{100}(0,8)(8,0)(16,8)
\end{picture} assuming that \ $y$ \ coordinate axis is horizontal.

There are two involutions on the parameter space of the standard versal deformation (see Table \ref{vd}) of this singularity: one of them takes the function \ $f_\lambda(x, y)$ \ to \ $f_\lambda(x,-y)$, and the other to \ $-f(-x,y)$. \ The first involution preserves the associated virtual morsifications (and hence both invariants \ $\mbox{\rm Card}$ \ and \ $\mbox{\rm Ind}$). The second involution preserves invariant \ $\mbox{\rm Card}$ \ but sends morsifications with \ $\mbox{\rm Ind}=I$ \ to these with \ $\mbox{\rm Ind}=-2-I$, in particular, it acts non-trivially on all components of \ $\R^{10} \setminus \Sigma$ \ (and on associated virtual components) with \ $\mbox{\rm Ind} \neq -1.$

\begin{figure}
\begin{picture}(40.00,35.00)
\put(0,5){\line(1,0){40}}
\put(0,5){\line(0,1){30}}
\put(0,35){\line(1,0){40}}
\put(40,5){\line(0,1){30}}
\bezier{300}(0,6)(20,20)(40,6)
\bezier{200}(0,21)(20,21)(0,34)
\bezier{200}(40,21)(20,21)(40,34)
\put(1,25){$-$} \put(18,6){$-$} \put(35,25){$-$}
\put(1,14){$+$}
\put(35,14){$+$}
\put(18,31){$+$}
\put(1,1){414538}
\put(36,1){0}
\end{picture} \qquad
\begin{picture}(40.00,37.00)
\put(0,5){\line(1,0){40}}
\put(0,5){\line(0,1){30}}
\put(0,35){\line(1,0){40}}
\put(40,5){\line(0,1){30}}
\bezier{300}(0,6)(20,20)(40,6)
\bezier{200}(0,21)(20,21)(0,34)
\bezier{200}(40,21)(20,21)(40,34)
\put(20,18){\circle{3}}
\put(1,25){$-$} \put(18,6){$-$} \put(35,25){$-$}
\put(1,14){$+$}
\put(35,14){$+$}
\put(18,31){$+$}
\put(1,1){308408}
\put(36,1){1}
\end{picture} \qquad
\begin{picture}(40.00,35.00)
\put(0,5){\line(1,0){40}}
\put(0,5){\line(0,1){30}}
\put(0,35){\line(1,0){40}}
\put(40,5){\line(0,1){30}}
\bezier{300}(0,6)(20,20)(40,6)
\bezier{200}(0,21)(20,21)(0,34)
\bezier{200}(40,21)(20,21)(40,34)
\put(17,18){\circle{3}}
\put(23,18){\circle{3}}
\put(1,25){$-$} \put(18,6){$-$} \put(35,25){$-$}
\put(1,14){$+$}
\put(35,14){$+$}
\put(18,31){$+$}
\put(1,1){179032} 
\put(36,1){2}
\end{picture} \qquad
\begin{picture}(40.00,35.00)
\put(0,5){\line(1,0){40}}
\put(0,5){\line(0,1){30}}
\put(0,35){\line(1,0){40}}
\put(40,5){\line(0,1){30}}
\bezier{300}(0,6)(20,20)(40,6)
\bezier{200}(0,21)(20,21)(0,34)
\bezier{200}(40,21)(20,21)(40,34)
\put(16,21){\circle{2.5}}
\put(24,21){\circle{2.5}}
\put(20,24){\circle{2.5}}
\put(1,25){$-$} \put(18,6){$-$} \put(35,25){$-$}
\put(1,14){$+$}
\put(35,14){$+$}
\put(18,31){$+$}
\put(1,1){21168}
\put(36,1){3}
\end{picture} \qquad
\begin{picture}(40.00,37.00)
\put(0,5){\line(1,0){40}}
\put(0,5){\line(0,1){30}}
\put(0,35){\line(1,0){40}}
\put(40,5){\line(0,1){30}}
\put(16,25){\circle{2.5}}
\put(24,25){\circle{2.5}}
\put(20,22){\circle{2.5}}
\bezier{300}(0,6)(20,23)(40,6)
\bezier{200}(0,21)(20,21)(0,34)
\bezier{200}(40,21)(20,21)(40,34)
\put(1,25){$-$} \put(18,6){$-$} \put(35,25){$-$}
\put(1,14){$+$}
\put(35,14){$+$}
\put(18,31){$+$}
\put(1,1){25416} 
\put(36,1){3}
\end{picture}
\qquad
\begin{picture}(40.00,37.00)
\put(0,5){\line(1,0){40}}
\put(0,5){\line(0,1){30}}
\put(0,35){\line(1,0){40}}
\put(40,5){\line(0,1){30}}
\bezier{300}(0,6)(20,20)(40,6)
\bezier{200}(0,21)(20,21)(0,34)
\bezier{200}(40,21)(20,21)(40,34)
\put(20,9){\circle{3}}
\put(1,25){$-$} \put(28,6){$-$} \put(35,25){$-$}
\put(1,14){$+$}
\put(35,14){$+$}
\put(18,31){$+$}
\put(1,1){18040}
\put(33,1){$-1$}
\end{picture} \qquad
\begin{picture}(40.00,35.00)
\put(0,5){\line(1,0){40}}
\put(0,5){\line(0,1){30}}
\put(0,35){\line(1,0){40}}
\put(40,5){\line(0,1){30}}
\bezier{200}(0,20)(20,20)(0,6)
\bezier{200}(40,21)(20,21)(40,34)
\bezier{200}(0,34)(10,23)(20,20)
\bezier{200}(20,20)(30,17)(40,6)
\put(1,25){$-$} \put(18,6){$-$} \put(35,25){$-$}
\put(1,14){$+$}
\put(35,14){$+$}
\put(18,31){$+$}
\put(1,1){336018}
\put(33,1){$-1$}
\end{picture} \qquad
\begin{picture}(40.00,35.00)
\put(0,5){\line(1,0){40}}
\put(0,5){\line(0,1){30}}
\put(0,35){\line(1,0){40}}
\put(40,5){\line(0,1){30}}
\bezier{200}(0,20)(20,20)(0,6)
\bezier{200}(40,21)(20,21)(40,34)
\bezier{200}(0,34)(10,23)(20,20)
\bezier{200}(20,20)(30,17)(40,6)
\put(22,23){\circle{3}}
\put(1,25){$-$} \put(18,6){$-$} \put(35,25){$-$}
\put(1,14){$+$}
\put(35,14){$+$}
\put(18,31){$+$}
\put(1,1){182596}
\put(36,1){0}
\end{picture} \qquad
\begin{picture}(40.00,37.00)
\put(0,5){\line(1,0){40}}
\put(0,5){\line(0,1){30}}
\put(0,35){\line(1,0){40}}
\put(40,5){\line(0,1){30}}
\bezier{200}(0,20)(20,20)(0,6)
\bezier{200}(40,21)(20,21)(40,34)
\bezier{200}(0,34)(10,23)(20,20)
\bezier{200}(20,20)(30,17)(40,6)
\put(16,25){\circle{3}}
\put(26,21){\circle{3}}
\put(1,25){$-$} \put(18,6){$-$} \put(35,25){$-$}
\put(1,14){$+$}
\put(35,14){$+$}
\put(18,31){$+$}
\put(1,1){103524}
\put(36,1){1}
\end{picture} \qquad
\begin{picture}(40.00,37.00)
\put(0,5){\line(1,0){40}}
\put(0,5){\line(0,1){30}}
\put(0,35){\line(1,0){40}}
\put(40,5){\line(0,1){30}}
\bezier{200}(0,20)(20,20)(0,6)
\bezier{200}(40,21)(20,21)(40,34)
\bezier{200}(0,34)(10,23)(20,20)
\bezier{200}(20,20)(30,17)(40,6)
\put(16,26){\circle{3}}
\put(26,26){\circle{3}}
\put(21,23){\circle{3}}
\put(1,25){$-$} \put(18,6){$-$} \put(35,25){$-$}
\put(1,14){$+$}
\put(35,14){$+$}
\put(18,31){$+$}
\put(1,1){40608}
\put(36,1){2}
\end{picture} \qquad
\begin{picture}(40.00,37.00)
\put(0,5){\line(1,0){40}}
\put(0,5){\line(0,1){30}}
\put(0,35){\line(1,0){40}}
\put(40,5){\line(0,1){30}}
\bezier{200}(0,20)(20,20)(0,6)
\bezier{200}(40,21)(22,20)(40,34)
\bezier{200}(0,34)(10,23)(20,20)
\bezier{200}(20,20)(30,17)(40,6)
\put(21,25){\circle{2}}
\put(23,22){\circle{2}}
\put(27,20){\circle{2}}
\put(30,19){\circle{2}}
\put(1,25){$-$} 
\put(18,6){$-$} 
\put(36,25){$-$}
\put(1,14){$+$}
\put(35,14){$+$}
\put(18,31){$+$}
\put(1,1){7200}
\put(36,1){3}
\end{picture} \qquad
\begin{picture}(40.00,35.00)
\put(0,5){\line(1,0){40}}
\put(0,5){\line(0,1){30}}
\put(0,35){\line(1,0){40}}
\put(40,5){\line(0,1){30}}
\bezier{200}(0,20)(20,20)(0,6)
\bezier{200}(40,21)(20,21)(40,34)
\bezier{200}(0,34)(10,23)(20,20)
\bezier{200}(20,20)(30,17)(40,6)
\put(22,23){\circle{3}}
\put(17,17){\circle{3}}
\put(1,25){$-$} 
\put(18,6){$-$} 
\put(35,25){$-$}
\put(1,14){$+$}
\put(35,14){$+$}
\put(18,31){$+$}
\put(1,1){12016}
\put(33,1){$-1$}
\end{picture} \qquad 
\begin{picture}(40.00,37.00)
\put(0,5){\line(1,0){40}}
\put(0,5){\line(0,1){30}}
\put(0,35){\line(1,0){40}}
\put(40,5){\line(0,1){30}}
\bezier{300}(0,6)(20,20)(40,6)
\bezier{300}(0,34)(20,20)(40,34)
\put(0,20){\line(1,0){40}}
\put(1,25){$-$} 
\put(18,6){$-$} 
\put(35,25){$-$}
\put(1,15){$+$}
\put(35,15){$+$}
\put(18,31){$+$}
\put(1,1){21040}
\put(33,1){$-1$}
\end{picture}
\caption{Perturbations of $J_{10}^3$}
\label{J103}
\end{figure}
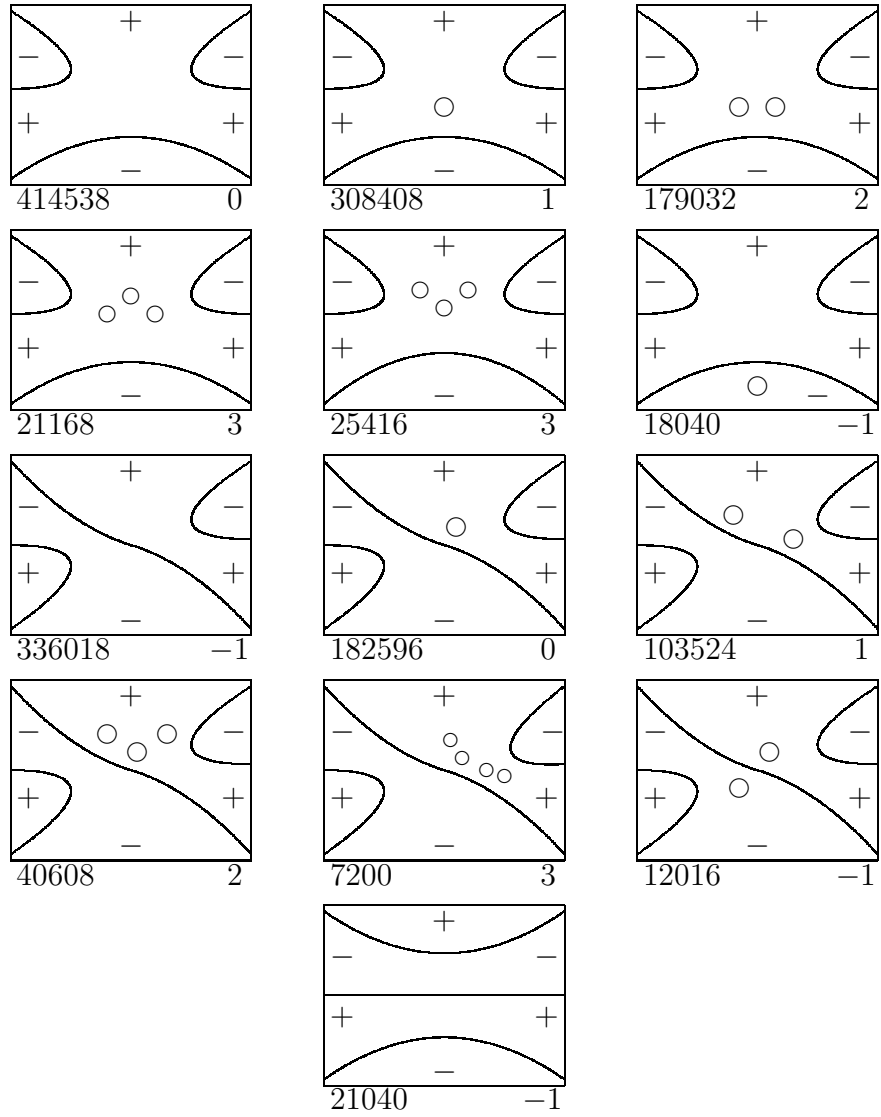

\begin{proposition}
\label{proJ3} 
1. Any singularity of type $J_{10}^3$ has exactly 2970134 distinct virtual morsifications divided into 23 virtual components.

2. There is an involution on the set of these virtual components that preserves invariant \ $\mbox{\rm Card}$ \ and maps any virtual component with \ $\mbox{\rm Ind} = I$ \ to a component with \ $\mbox{\rm Ind} = -I -2,$ \ in particular acts non-trivially on all virtual components with $\mbox{\rm Ind} \neq -1$. There are 13 orbits of this action.

3. Topological types of some realizations of representatives of all these orbits are shown in Fig.~\ref{J103}; the values of invariants \ $\mbox{\rm Card}$ \ and \ $\mbox{\rm Ind}$ \ are indicated there by left- and right-hand subscripts. The realizations of virtual components invariant under our involution are shown in all the pictures with \ $\mbox{\rm Ind} = -1$ \ except for the one with \ $\mbox{\rm Card} = 18040$. The pictures of realizations of remaining ten virtual components can be obtained from pictures of their mates by reflecting in a horizontal mirror and changing all signs. 

4. There are additional ten components of the complement of the discriminant of \ $J_{10}^3$ \ singularity, the realizations of which can be obtained by the reflection in the vertical mirror $($i.e. by the change \ $f_\lambda(x,y) \mapsto f_\lambda(x,-y))$ \ from pictures of Fig.~\ref{J103} with \ $\mbox{\rm Card}$ \ equal to 182596, 103524, 40608, 7200, 336018, 12016, and from the pictures of the mates of first four of them by the involution described in statement 2.  
\end{proposition}

\begin{remark} \rm
All 33 components of the complement of the discriminant variety mentioned in this proposition are separated also by the topological type invariant (see \S \ref{invariants}), except only for two pairs of components  with the values of \ $\mbox{\rm Card}$ \ equal to 21168 and 25416, which are separated by Bezout's theorem considerations.
\end{remark}

\begin{conjecture}
\label{conjJ3}
Any real function singularity of type \ $J_{10}^3$ \ has exactly 33 different components of the complement of discriminant variety in the parameter space of any its versal deformation.
\end{conjecture}

\unitlength=1.0mm
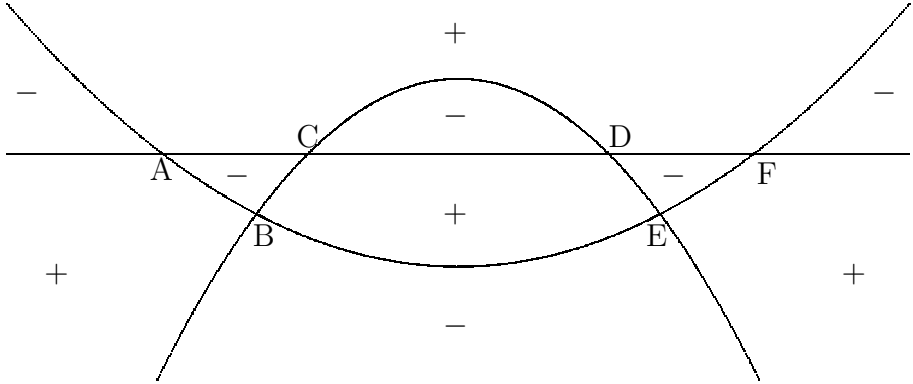
\begin{figure}
\begin{picture}(120,50)
\put(0,30){\line(1,0){120}}
\bezier{500}(0,50)(60,-20)(120,50)
\bezier{500}(20,0)(60,80)(100,0)
\put(58,21){$+$}
\put(29,26){$-$}
\put(87,26){$-$}
\put(58,34){$-$}
\put(19,26.5){A}
\put(99.7,26){F}
\put(32.7,17.8){B}
\put(85,17.8){E}
\put(38.5,31){C}
\put(80,31){D}
\put(58,6){$-$}
\put(58,45){$+$}
\put(1,37){$-$}
\put(115,37){$-$}
\put(5,13){$+$}
\put(111,13){$+$}
\end{picture}
\caption{First perturbation for $J_{10}^3$}
\label{J103z}
\end{figure}

\noindent
{\it Proof of Proposition \ref{proJ3}.}
Let us take a perturbation of \ $f_0$, whose zero set is shown in Fig.~\ref{J103z}. Let \ $\tilde f$ \ be its very small further perturbation, which is a strictly Morse function with critical values at saddlepoints ordered as \ $\tilde f(A)<\tilde f(B)<\tilde f(C)<\tilde f(D)<\tilde f(E)<\tilde f(F)$.

The first eight elements of the standard scale of \ $\tilde f$ \ realize pictures of Fig.~\ref{J103} with left subscripts equal to 414538, 308408, 179032, 21168, 179032, 103524, 182596, 336018. Then we take another small perturbation \ $\tilde f_1$ \ of the same picture of Fig.~\ref{J103z} such that critical values at saddlepoints are ordered as \ $\tilde f_1(E)<\tilde f_1(B)<\tilde f_1(C)<\tilde f_1(D)<\tilde f_1(A)<\tilde f_1(F)$. The function \ $\tilde f_1 - c$ \ of its standard scale, which has exactly seven  negative critical values, realizes the picture  with left subscript  $18040$. 
Let $\tilde f_2$ be one more strictly Morse perturbation of the function of Fig.~\ref{J103z}, chosen so that  $\tilde f_2(A)<\tilde f_2(C)<\tilde f_2(F)<\tilde f_2(D)<\tilde f_2(B)<\tilde f_2(E)$. The function $\tilde f_2 - c$ of its standard scale, which has exactly seven critical values, realizes the picture with left subscript $21040$.

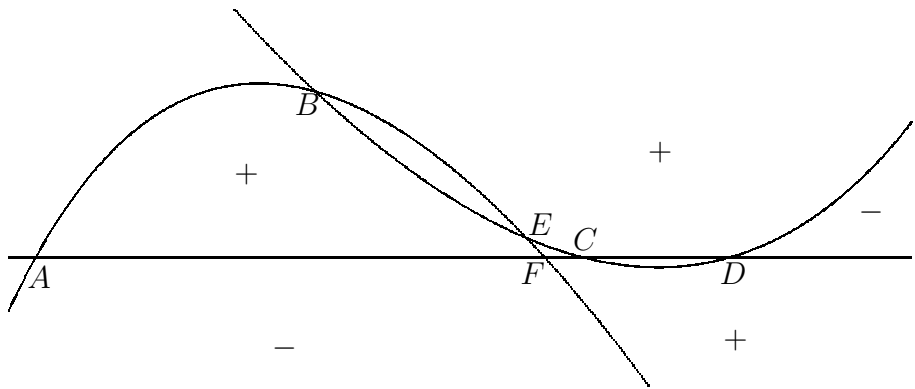
\begin{figure}
\begin{picture}(120,50)
\put(0,17){\line(1,0){120}}
\bezier{500}(30,50)(85,-10)(120,35)
\bezier{500}(0,10)(30,75)(85,0)
\put(2.5,13){$A$}
\put(38,36){$B$}
\put(69,20){$E$}
\put(68,13.5){$F$}
\put(75,18){$C$}
\put(94.5,13.5){$D$}
\put(30,27){$+$}
\put(85,30){$+$}
\put(35,4){$-$}
\put(95,5){$+$}
\put(113,22){$-$}
\end{picture}
\caption{Second perturbation for $J_{10}^3$}
\label{J103a}
\end{figure}

Next, we take another starting perturbation of \ $f_0$, whose set of zeros is as shown in  Fig.~\ref{J103a}. Let \ $\tilde f_3$ \ be its slight further perturbation such that \ $\tilde f_3(A)<\tilde f_3(B)<\tilde f_3(C)<\tilde f_3(D)<\tilde f_3(E)<\tilde f_3(F)$. The function \ $\tilde f_3 - c$ \ of its standard scale, which has exactly three (respectively, four) negative critical values, realizes the picture of Fig.~\ref{J103} with subscript $25416$  (respectively,  $40608$). Let \ $\tilde f_4$ \ be yet another Morse perturbation of Fig.~\ref{J103a} such that \ $\tilde f_4(E) < \tilde f_4(A) < \tilde f_4(B) < \tilde f_4(F) < \tilde f_4(D) < \tilde f_4(C) $. Then the function \ $\tilde f_4 - c $ \ with exactly seven negative critical values  realizes the picture  with left subscript $12016$. 

The virtual morsifications associated with all twelve constructed perturbations can be found by the method of \cite{GZ}, \cite{AC}. Substituting them into our reduced program, we find that these left subscripts are indeed equal to the corresponding values of invariant \ $\mbox{\rm Card}$.

The involution on the set of virtual components promised in statement 2 of Proposition \ref{proJ3} reflects the involution \ $f_\lambda(x,y) \leftrightarrow -f_\lambda(-x,y)$: \ it reverses the order of the critical values (and the corresponding vanishing cycles), turns minima to maxima and vice versa, and replaces the number of negative critical values by that of positive ones. Applying it to all twelve virtual morsifications realized above, in nine cases we arrive at a different virtual component. Namely, in eight cases with \ $\mbox{\rm Ind} \neq -1$ \ this involution changes invariant \ $\mbox{\rm Ind}$, \ and the virtual component of the morsification with \ $\mbox{\rm Card} = 18040$ \ does not contain virtual morsifications with only three negative critical values (as evidenced by our reduced program on a special request), while the result of our involution applied to this morsification has exactly three negative critical values.

The main program tells us that the total number of virtual morsifications of \ $f_0$ \ is  2970134. Summing up the cardinalities of all 21 virtual components constructed above, we see that there are exactly 14400 virtual morsifications not covered by them.

One of these missing virtual morsifications is characterized by the intersection matrix 
\begin{equation} \left|
\begin{array}{cccccccccc}
$-2$  &  0  &  0   & 0 &   0  &  1  &  0  &  0  &  0  &  0 \\
 0  & $-2$   &  0   & 0  &  0  &  0  &  1 &   0  &  1  &  0 \\ 
 0  &  0  & $-2$   & 0 &  $-2$  &  1  &  0  &  1  &  0   & 1 \\
0  &  0  &  0   & $-2$  &  0 &   0  &  1 &   0  &  0  &  1 \\
0  &  0  & $-2$  &  0 &  $-2$  &  1  &  0  &  1   & 0   & 1 \\
1  &  0  &  1   & 0   & 1  & $-2$  &  0  &  0  &  0  &  0  \\
0   & 1  &  0   & 1  &  0   & 0  & $-2$  &  0  &  0   & 0  \\
0 &   0  &  1  &  0   & 1  &  0  &  0  & $-2$  & 0  &  0  \\
0   & 1  &  0 &   0  &  0  &  0   & 0  &  0 &  $-2$  &  0   \\
0  &  0  &  1  &  1  &  1  &  0  &  0   & 0  &  0  & $-2$
\end{array}
\right| \ ,\label{matJ10}
\end{equation}
the string  (0, 0, 1, 0, 0, 1, 1, 1, 1, 1) of Morse indices, and the number of negative critical values equal to five. Our reduced program shows that its virtual component consists of 7200 virtual morsifications; the value of invariant \ $\mbox{\rm Ind}$ \ on it is equal to $3 \neq -1$, therefore our involution moves it to another virtual component, also with \ $\mbox{\rm Card}=7200$. So, these 23 virtual components cover all virtual morsifications of our singularity. This finishes a proof of statements 1 and 2 of Proposition \ref{proJ3}.

In order to realize two last virtual morsifications and  the picture of Fig.~\ref{J103} marked with subscript $7200$,  we first simplify our singularity \ $J_{10}^3$ \ to \ $E_8$\ (by adding the term \ $\varepsilon y^5$ \ to the normal form of Table \ref{t4}) and further decompose the obtained singularity of type \ $E_8$ \ into four minima and four saddlepoints, as shown in the fifth picture of Fig.~7 of \cite{VLoo}, see also Fig.~12 there and page 17 of \cite{AC}. 

The corresponding virtual morsification can belong only to a virtual component with \ $\mbox{\rm Card}=7200$. Indeed, invariant \ $\mbox{\rm Ind}$ takes value 3 on it. For both two competing virtual components with this value, i.e., having \ $\mbox{\rm Card}$ \ equal to 21168 or 25416, our restricted program testifies (upon a special request) that they do not contain virtual morsifications whose second, third, fourth and fifth critical points (in ascending order of critical values) are minima, as they are for the morsification just constructed.

Finally, statement 4 of Proposition \ref{proJ3} follows directly from these realizations. \hfill $\Box$

\section{$P_8^1$}
This class is represented by non-degenerate homogeneous cubic polynomials in \ $\R^3$ \ such that corresponding curves in \ $\RP^2$ \ consist of only one component. 

The involution 
\begin{equation}
\label{invol}
f_\lambda (x,y,z) \leftrightarrow -f_\lambda(-x,-y,-z)
\end{equation} 
acts on the parameter space $\R^8$ of the versal deformation of Table \ref{vd} of such a singularity; it sends the discriminant to itself and hence permutes somehow the components of its complement.

\begin{table}
\caption{Statistics for $P_8^1$ singularity}
\label{p8stat}
\begin{tabular}{|c|ccccccc|c|}
\hline
No. & 1 & 2 & 3 & 4 & 5 & 6 & 7 \\
\hline
Ind & $ -3$ & $ -3$ & $-2$ & $-1$ & 0 & 1 & 1 \\
 Card & 210 & 210 & 1370 & 1992 & 2465 & 128 & 128 \\
\hline 
\end{tabular}
\end{table}

\begin{proposition}
\label{proP1}
1. Any  singularity of type  $P_8^1$ has exactly 
6503 different virtual morsifications, which are divided into seven virtual components.

2. Invariants \ $\mbox{\rm Ind}$ \ and \ $\mbox{\rm Card}$ \ of these virtual components are shown in  Table \ref{p8stat}.

3. All these virtual components are associated with $($i.e. contain virtual morsifications of real morsifications from$)$  some real components of the complement of the discriminant of any versal deformation of a singularity of class $P_8^1$.

4. The closure in $\RP^3$ of the zero set of any non-discriminant perturbation $f_\lambda$ of our singularity with the value of invariant \ $\mbox{\rm Ind}$ \ equal to 0 $($respectively, $-1$, $-2$, $-3$, and 1$)$ is homeomorphic to the projective plane $($respectively, the projective plane with one, two, three handles, and the disjoint union of the projective plane and the sphere$)$, so  that they realize all topological types of non-singular cubic surfaces in $\RP^3$, see \cite{schlaf}, \cite{klein}, \cite{zeuten}.
\end{proposition}

\begin{remark} \rm
\label{remP}
This proposition agrees with the result of \cite{FK}, which states  that there are exactly five rigid isotopy classes of real affine smooth cubic surfaces transversal to the infinitely distant plane along a single curve, and these classes represent all five classes of smooth cubic surfaces in $\RP^3$ indicated in \cite{schlaf}, \cite{klein}, \cite{zeuten}. We have more than five virtual components because we consider functions and not their zero sets only, so that e.g. morsifications related by involution (\ref{invol}) can represent different components, although their zero sets (and also closures of these sets) are equivalent. 

Moreover, dealing with deformations of singularities we consider functions with a fixed (or at least only slightly varying) principal homogeneous parts, so a priori there can be several real components of $\R^8 \setminus \Sigma$ associated with a virtual component. Say, they can be reducible to one another by affine transformations of $\R^3$ from the group  (isomorphic to $S(3)$) of symmetries of the principal homogeneous part of $f_0$.
 For this reason, in this case we are not sure that the conjecture about the number of  components of complements of discriminants formulated at the end  of Theorem \ref{mtreal} is correct (unlike the cases $X_9$ and $J_{10}$).
\end{remark}

\subsection{Proof of statements 1 and 2 of Proposition \ref{proP1}} 
Take the simplest representative $f_0= x^3 + y^3 + z^3$ of class $P_8^1$ and consider its  morsification $\tilde f \equiv f_0 - \varepsilon_1 x - \varepsilon_2 y - \varepsilon_3 z$, where all coefficients $\varepsilon_i$ are positive, distinct, and sufficiently close to one another. The intersection form of vanishing cycles for $\tilde f$ can be calculated by the main theorem of \cite{gab}, see \S II.2.3 in \cite{AGLV1}; the remaining ingredients of the corresponding virtual morsification follow easily from it.
Our main program with these initial data says that 
\begin{enumerate}
\item 
Singularity $P_8^1$ has exactly 6503 virtual morsifications,
\item 
invariant \ $\mbox{\rm Ind}$ \ takes on them values from $-3$ to $+1$, namely, there are respectively 420, 1370, 1992, 2465 and 256 distinct virtual morsifications with  values of  \ $\mbox{\rm Ind}$ \ equal to $-3$, $-2$, $-1$, 0 and $1$.
\end{enumerate}

Next, we take the standard scale of functions $\tilde f - c$ \ for this morsification $\tilde f$ and apply the reduced program to each of its nine elements. The pairs of values of invariants \ $(\mbox{\rm Card}, \mbox{\rm Ind})$ \ for these elements include $(128,1), $ $(2465,0), $ $(1992, -1),$ and $ (1370,-2)$ (but morsifications with \ $\mbox{\rm Ind}= -3$ \ do not appear in this scale). Comparing this with the results of the main program, we see that each of the values 0, $-1$ and $-2$ of  invariant $\mbox{\rm Ind}$ is represented by a single virtual component. The second and eighth elements of this scale are related to each other via the involution (\ref{invol}) 
and have equal values $(128,1)$ of invariants \ $(\mbox{\rm Card}, \mbox{\rm Ind})$. \ Nevertheless, the  components of $\R^8 \setminus \Sigma$ containing these elements are different. This follows from the topology of corresponding sets of lover values, $f_\lambda^{-1}((-\infty,0])$: for one of these elements  this set is homotopy equivalent to $S^2$, and for the other to two points.  

Our reduced program also proves the difference of  corresponding two virtual components. Namely, it detects (upon a special request) that one of these virtual components does not contain virtual morsifications with seven negative critical values, and the other component does not contain virtual morsifications with seven positive ones (while their starting virtual morsifications have respectively seven positive/negative critical values). Since $128+128 = 256$, we have found and realized all virtual components with $\mbox{\rm Ind}=1$.

Then we ask the main program to find and print out a virtual morsification with $\mbox{\rm Ind}= -3$ and only real critical points. A possible answer is presented by matrix (\ref{matp81}), string of Morse indices equal to (even, odd, odd, odd, odd, even, even, even), and the number of negative critical values equal to five.
\begin{equation} \left|
\begin{array}{cccccccc}
    $-2$  & 1 &  1 &  1 &  1  &  $-1$  &  $-1$  &  $-1$ \\
   1  &  $-2$  &  0   & 0 &   0 &  1  &  0  &  0 \\
   1  &  0   & $-2$  &  0  &  0 &   0  &  0 &  1 \\
   1  &  0   & 0   & $-2$  &  0   & 0  & 1  &  0 \\
   1  &  0  &  0  &  0  &  $-2$  & 1  & 1  & 1 \\
    $-1$ &  1  &  0  &  0 &  1  &  $-2$  &  0  &  0 \\
    $-1$  &  0  &  0 &  1 &  1  &  0   & $-2$  &  0 \\
    $-1$  &  0  & 1  &  0  & 1   & 0  &  0  &  $-2$
\end{array}
\right| \label{matp81}
\end{equation}
Substituting it as the initial data into the reduced program, we get that the  corresponding virtual component consists of 210 different virtual morsifications. 

Further, we take the ``negative'' of this virtual morsification (obtained from it by the change of virtual data modeling the involution (\ref{invol})) and run the reduced program with these initial data. It again finds some 210 virtual morsifications; it remains to prove that this list  does not coincide with the previous one. To do it, in both cases I asked the reduced program to print out all obtained virtual morsifications {\em not all of whose critical points are real}. It turned out that each of these two lists consists of only 42 elements, and it is easy to check by hand that they do no intersect. Since (according to the main program) all virtual components with $\mbox{\rm Ind} = -3$ contain exactly 420 virtual morsifications, we conclude that we have found all of them. \hfill $\Box$ 

\begin{problem} \rm
Construct an algebraic or combinatorial invariant of virtual morsifications separating (in particular) the virtual components discussed in the last part of the previous proof.
\end{problem}

\subsection{Realization of virtual components} Some five out of seven virtual components from Proposition \ref{proP1} are already realized by the elements of the standard scale of morsification $\tilde f$. It remains to find a perturbation $f_\lambda$ of our singularity $f_0 \in P_8^1$, such that $\mbox{\rm Ind}(f_\lambda)= -3$. The function related with it by involution (\ref{invol}) will then represent a different component of the complement of the discriminant.

Let the original singularity $f_0 \in P_8^1$ again have the form $x^3 + y^3 + z^3$, and let $l:\R^3 \to \R$ be a linear form whose projectivized zero set in $\RP^2$  is tangent to that of $f_0$ at one of  its three inflection points. Then all functions of the form $f_0 + \tau l^2$, $\tau \neq 0$, belong to the singularity class $E_6$. Namely, if $\tau >0$, then such a function has a critical point equivalent to the polynomial $\tilde x^3-\tilde y^4 + \tilde z^2 \in -E_6$ with critical value 0, and a critical point of type $A_2$ with positive critical value. 
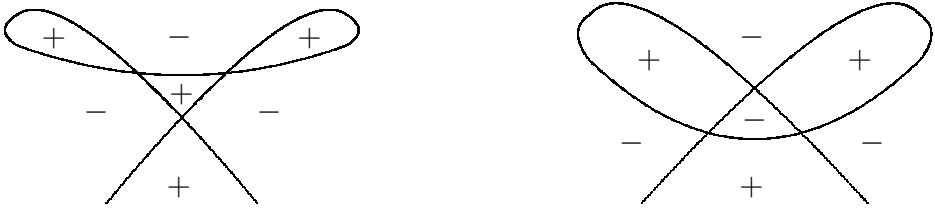
\begin{figure}
\begin{center}
{
\begin{picture}(50,26)
\bezier{300}(15,1)(40,31)(47,26)
\bezier{70}(47,26)(50,24)(47,22)
\bezier{250}(47,22)(25,14)(3,22)
\bezier{70}(3,22)(0,24)(3,26)
\bezier{300}(3,26)(10,31)(35,1)
\put(25.00,15.50){\makebox(0,0)[cc]{$+$}}
\put(42.00,23.00){\makebox(0,0)[cc]{$+$}}
\put(8.00,23.00){\makebox(0,0)[cc]{$+$}}
\put(23,2){$+$}
\put(23,22){$-$}
\put(12,12){$-$}
\put(35,12){$-$}
\end{picture}
\qquad \qquad \qquad
\begin{picture}(50,26)
\bezier{300}(10,1)(40,34)(47,26)
\bezier{70}(47,26)(50,24)(47,20)
\bezier{300}(47,20)(25,-1)(3,20)
\bezier{70}(3,20)(0,24)(3,26)
\bezier{300}(3,26)(10,34)(40,1)
\put(25.00,12.00){\makebox(0,0)[cc]{$-$}}
\put(39.00,20.00){\makebox(0,0)[cc]{$+$}}
\put(11.00,20.00){\makebox(0,0)[cc]{$+$}}
\put(23,2){$+$}
\put(23,22){$-$}
\put(7,8){$-$}
\put(39,8){$-$}
\end{picture}
}
\end{center}
\caption{Perturbations for $-E_6$}
\label{e6proof}
\end{figure}

According to \cite{GZ} or \cite{AC}, this critical point  of type $-E_6$ can be additionally perturbed so that it splits into three points with Morse index 1 and critical value 0, and three points with Morse index 2 and positive critical values; see Fig.~\ref{e6proof} (left) for the corresponding two-dimensional picture. Subtracting from the obtained function a positive constant smaller than the latter three critical values, we obtain a function with $\mbox{\rm Ind} = -3.$
\medskip

Statement 4 of Proposition \ref{proP1} follows immediately from the realizations. \hfill $\Box$

\section{$P_8^2$}
\label{P82}

This class is represented by any non-degenerate homogeneous polynomials $f_0$ of degree 3 in $\R^3$, whose projectivized zero sets in $\RP^2$ consist of one oval and one component isotopic to the subspace $\RP^1 \subset \RP^2$. We can and will assume that our  versal deformation of $f_0$ has the form indicated in Table \ref{vd}.

\subsection{A homological  invariant}
\label{hominv}
We will use the following simplification of the ``topological type'' invariant from \S \ref{invariants}. Let $D$ be a ball centered at the origin in the parameter space $\R^8$ of our versal deformation. Let $B$ be a sufficiently large (compared to $D$) ball in $\R^3$. Then for any $\lambda \in D$ the intersection of the set of lower values $W_\lambda \equiv  f_\lambda^{-1}((-\infty, 0]) $ and  sphere $ \partial B$ is homeomorphic to the disjoint union of a disc and an annulus, in particular, its reduced homology group $\tilde H_i(W_\lambda \cap \partial B)$ is equal to $\Z$ for $i$ equal to 0 or 1 and is trivial for all other $i$. 

\begin{definition} \rm 
{\it Homology index} of  function $f_\lambda$, $\lambda \in D$ (and of entire component of $\R^8 \setminus \Sigma$ containing $\lambda$) is the sequence of  Betti numbers of relative homology groups $\tilde H_j(W_\lambda \cap B, W_\lambda \cap \partial B)$, $j=0, 1, 2$, the last  two of which can additionally be marked by bars above them if the boundary operator from the corresponding group to $\tilde H_{j-1}(W_\lambda \cap \partial B)$ is non-trivial. 
\end{definition} 

\begin{lemma} 
\label{lemev}
1. Invariant \ $\mbox{\rm Ind}$ \ is equal to the alternating sum of three elements of the homology index. 

2. The bar stands over the second $($respectively, the third$)$ element of the homology index of a real component if and only if does not stand over the third $($respectively, second$)$ element of the homology index of the component related with the initial component by involution $($\ref{invol}$)$. In particular, if a component of the complement of the discriminant is invariant under this involution, then exactly one of its elements is equipped with a bar. 

3. The third element of the homology index of  a component of $\R^8 \setminus \Sigma$  or of the component related to it by involution $($\ref{invol}$)$ is  barred if and only if for any $\lambda$ from this component the ``oval''  part of the intersection of the projective closure of corresponding zero level set $f^{-1}_\lambda(0)$ with ``infinity'' plane $\RP^3 \setminus \R^3$ is zero-homologous in this projective closure. 
\end{lemma}

\noindent
{\it Proof.} All these facts are elementary. \hfill $\Box$

\begin{table}
\caption{Statistics for $P_8^2$ singularity}
\label{p82stat}
\noindent
{\footnotesize
\begin{tabular}{|c|cccc|ccc|ccc|cc|ccc|c|}
\hline
No. & 1 & 2 & 3 & 4 & 5 & 6 & 7 & 8 & 9 & 10 & 11 & 12 & 13 & 14 & 15  \\
\hline
Ind & \multicolumn{4}{c|}{$-3$} & \multicolumn{3}{c|}{$-2$} & \multicolumn{3}{c|}{$-1$} & \multicolumn{2}{c|}{0} & \multicolumn{3}{c|}{1}  \\
\hline
Card & 258 & 156 & 60 & 60 & 1216 & 336 & 336 & 1318 & 844 & 844 & 1648 & 1648 & 262 & 94 & 94  \\
HI & $0 \bar 3 0$ & $0 \bar 3 0$ & $030$ & $0\bar 4 \bar 1$ & $0 \bar 2 0$ & $020$ & $0\bar 3 \bar 1$ & $0 \bar 1 0$ & 010 & $0 \bar 2 \bar 1$ & $000$ & $0\bar 1 \bar 1$ & $00\bar 1 $ & $100$ & $0 \bar 1 \bar 2$  \\ 
\hline 
\end{tabular}
}
\end{table}

\begin{proposition}
\label{proP2}
1. Each singularity $f_0$ of class $P_8^2$ has exactly 9174 distinct virtual morsifications. 

2. These virtual morsifications are distributed over exactly fifteen virtual components. The values of  invariants \ $\mbox{\rm Ind}$ \ and \ $ \mbox{\rm Card}$ \ of these components are shown in the second and third rows of Table \ref{p82stat}. 

3.  The closure in $\RP^3$ of the set of zeros of any non-discriminant perturbation $f_\lambda$ of our singularity with value of \ $\mbox{\rm Ind}$ \ equal to 0 $($respectively, $-1$, $-2$, $-3$, and 1$)$ is homeomorphic to the projective plane $($respectively, projective plane with one, two, three handles, and the disjoint union of the projective plane and the sphere$)$.  

4. All pairs of different virtual components of $P_8^2$ singularity with the same values of invariant \ $\mbox{\rm Card}$ \ are realized by morsifications connected by involution $($\ref{invol}$)$.
\end{proposition}

A proof of this proposition takes the rest of this section. 

\begin{remark} \rm
This proposition agrees with one more result of the work \cite{FK}, which enumerates  rigid isotopy classes of smooth affine cubic surfaces in $\R^3$  intersecting the ``infinitely distant'' plane transversally along two connected curves. As in the case of $P_8^1$, we have more classes than \cite{FK} because we consider functions and not just their zero level sets, so in some cases perturbations $f_\lambda(x,y,z)$ and $-f_\lambda(-x,-y,-z)$ belong to different components of $\R^8 \setminus \Sigma$ but have equivalent zero level sets.  Namely, this happens for all five pairs of columns of Table \ref{p82stat} that have equal values of both invariants $\mbox{\rm Ind}$ \ and \ $\mbox{\rm Card}$.

And again, I do not know whether a virtual component of a $P_8^2$ singularity can be associated with several components of $\R^8 \setminus \Sigma$ moved into one another by the action of the symmetry group (isomorphic to $S(3)$) of singularities $P_8$ and their versal deformations.
\end{remark}

\begin{problem}
To fix the correspondence between two virtual components from Table \ref{p82stat} with the values of \ $\mbox{\rm Card}$ \ equal to 258 and 156 on one hand, and the isotopy  classes marked in Fig.~ 1 of \cite{FK} with the pairs of numbers $(15, 12)$ and $(11,16)$, on the other. $($The correspondence between all other elements of these lists, associating with any virtual morsification the isotopy class of surfaces representing it, is trivial$)$.
\end{problem}

\noindent
{\it Proof of Proposition \ref{proP2}.} 
\label{proofP2}
Our starting point is the morsification $\tilde f$ of a  singularity of class $P_8^2$ found in \S 7 of \cite{VMS}. Namely, we take the function $f_0 \in P_8^2$ with Newton--Weierstrass normal form equal to $x^3 - x z^2 + y^2 z$. Similarly to the previous section, we consider a family of its small perturbations of the form $f_\tau = f_0 + \tau z^2,$ which, for $\tau \neq 0$, have critical points of the class $+E_6$ or $-E_6$ at the origin. In particular, for $\tau>0$ this critical point is of type $-E_6$ and  is equal to $\tilde x^3 - y^4 + \tilde z^2$ in some local coordinates. This function $f_\tau$, $\tau \neq 0$, also has two Morse critical points with Morse indices 1 and 2 and equal critical values, the sign of which is equal to the sign of $\tau$. Fix such a perturbation $f_{\tau}$ with some small value $\tau>0$ and take an additional tiny perturbation of $f_\tau$ that moves the critical point of type $-E_6$ as shown in  Fig.~\ref{e6proof} (left), i.e. splits it into three Morse points with Morse index 1 and critical value 0 and three points with Morse index 2 and  positive critical values. In addition, we require that critical values of two Morse critical points arising near the Morse points of $f_\tau$ become different: namely, the value at the point with Morse index 1 becomes greater  than that at the point with Morse index 2. Finally, we perturb the obtained function once more in such a way that all critical values at three points of index 1 corresponding to crossing points in Fig.~\ref{e6proof} (left) become different. The resulting function is the desired morsification $\tilde f$.  

It was proved in \cite{VMS} that the intersection form of vanishing cycles of this morsification (ordered by the increase of corresponding critical values)
is given by matrix (\ref{matp8}); Morse indices of the corresponding critical points are equal to 1, 1, 1, 2, 2, 2, 2, and 1, respectively.
\begin{equation} \left|
\begin{array}{cccccccc}
$-2$ \ & 0 & 0 & 1 & 0 & 1 & 0 & 1 \\
0 & $-2$ \ & 0 & 0 & 1 & 1 & 0 & 1 \\
0 & 0 & $-2$ \ & 0 & 0 & 1 & 1 & 1 \\
1 & 0 & 0 & $-2$ \ & 0 & 0 & 0 & 0 \\
0 & 1 & 0 & 0 & $-2$ \ & 0 & 0 & 0 \\
1 & 1 & 1 & 0 & 0 & $-2$ \ & 0 & $-2$ \\
0 & 0 & 1 & 0 & 0 & 0 & $-2$ \ & 0 \\
1 & 1 & 1 & 0 & 0 & $-2$ & 0 & $-2$ \
\end{array}
\right| \label{matp8}
\end{equation}

We substitute this data into our main program, and it tells us the following facts:
\begin{lemma}
\label{lem77}
1. The total number of virtual morsifications of $P_8^2$ singularity is 9174.

2. Invariant \ $\mbox{\rm Ind}$ \ takes values from $-3$ to $+1$ on these virtual morsifications.

3. Numbers of virtual morsifications with the value of \ $\mbox{\rm Ind}$ \ equal to $-3,$ $ -2,$ $ -1,$ 0, and 1 are respectively $534,$ $1888,$ $3006,$ $ 3296,$ and $450$. \hfill $\Box$
\end{lemma}

Then we apply the reduced program, starting from nine functions $\tilde f -c$ of the standard scale of $\tilde f$. We get that values of two main invariants on these functions are equal to  $(1648,0), $ $(844,-1),$ $(336,-2),$ $(60,-3),$ $(336,-2),$ $(844,-1),$ $(1648,0),$ $(262,1)$, (1648,0). Although the first and the last elements of this scale have the same pair of invariants $(1648,0),$ corresponding real components of $\R^8 \setminus \Sigma$ are separated by their homology indices of \S \ref{hominv}: the set of lower values is homotopy equivalent to a pair of points in the first case and to $S^2$ in the second. Also, our reduced program detects that the virtual component of the first (respectively, second) of these two virtual morsifications does not contain virtual morsifications with eight negative (respectively, eight positive) critical values, while for the other virtual component  such a morsification is the starting point of its calculation. Therefore, these two virtual components also are different.

Similarly, the virtual morsifications of the second (respectively, third, fourth) element of this scale and of its negative (i.e., the function obtained from it by involution (\ref{invol})) belong to different virtual components. Indeed, our reduced program checks that the virtual component of this element  does not contain virtual morsifications with seven (respectively, six, five) negative critical values, and the virtual component of the negative of this element does not contain virtual morsifications with seven (respectively, six, five) positive critical values.

Further, constructing morsification $\tilde f$  we could change the order of the 7-th and 8-th critical values. Indeed, both corresponding critical points are obtained by perturbing points with the same positive critical value, and the intersection index of the corresponding vanishing cycles is 0. This gives us the same matrix as (\ref{matp8}) with only the last two columns and rows permuted, and Morse indices of the last two critical points become equal to (1, 2) and not (2, 1). Then we substitute this data into the reduced program (starting from the element of the corresponding standard scale, which has seven negative critical values) and get a new virtual component with \ $\mbox{\rm Card} =1318$ \ and \ $\mbox{\rm Ind} =-1$. 

Also, constructing the morsification $\tilde f$ (see page \pageref{proofP2}) we could split the  singularity of type $-E_6$ in accordance with the right-hand picture of Fig.~\ref{e6proof}, and not the left-hand one, which was previously used. Let $\tilde f_2$ be the resulting function. To find virtual morsification data of this function or at least from the same virtual component, we can take the function $\tilde f - c $ from the standard scale of $\tilde f$ having exactly six negative critical values, and ask the reduced program to find a virtual morsification in this virtual component, the lowest critical value of which is reached at a point with even Morse index. A virtual morsification obtained in this way is characterized by intersection matrix (\ref{matp8-3}), string of Morse indices (even, odd, odd, odd, even, even, even, odd), and only one negative critical value.
\begin{equation} \left|
\begin{array}{cccccccc}
$-2$ & 1 & 1 & 1 & $-1$ & $-1$ & $-1$ & 1 \\
1 & $-2$ & 0 & 0 & 1 & 1 & 0 & 0 \\
1 & 0 & $-2$ & 0 & 1 & 0 & 1 & 0 \\
1 & 0 & 0 & $-2$ & 0 & 1 & 1 & 0 \\
$-1$ & 1 & 1 & 0 & $-2$ & 0 & 0 & 0 \\
$-1$ & 1 & 0 & 1 & 0 & $-2$ & 0 & 0 \\
$-1$ & 0 & 1 & 1 & 0 & 0 & $-2$ & 0 \\
1 & 0 & 0 & 0 & 0 & 0 & 0 & $-2$ 
\end{array}
\right| \label{matp8-3}
\end{equation} 
Applying the reduced program with this initial data, we get that the corresponding virtual component consists of 94 virtual morsifications. This program also tells us that none of these virtual morsifications has seven negative critical values, therefore the negative of our virtual morsification does not belong to the same virtual component (and represents another virtual component, also of 94 elements). 

So, we have found (and realized) twelve distinct virtual components with the values of the pair of invariants \ $(\mbox{\rm Card}, \mbox{\rm Ind})$ \ equal to $(94,1), (94, 1), (262,1),$ $ (1648,0), $ $(1648,0), (844,-1),$ $ (844,-1),$ $ (1318,-1), $ $(336,-2), $ $(336,-2),$ $ (60,-3),$ $ ( 60,-3)$. Comparing this with Lemma \ref{lem77}, we see that the found virtual components cover all virtual morsifications with  \ $\mbox{\rm Ind}$ \ equal to $1, 0$ and $-1$, but  1216 virtual morsifications with \ $\mbox{\rm Ind}=-2$ and 414 virtual morsifications with \ $\mbox{\rm Ind}=-3$ remain uncovered.

Some two of virtual morsifications calculated by the main program have intersection matrices (\ref{matp8-2}); the corresponding Morse indices in the first case are (odd, odd, even, odd, odd, even, even, even) and in the second case (odd, odd,odd, even, even, even, odd, even).
\begin{equation} \left|
\begin{array}{cccccccc}
$-2$ \ & 0 & 1 & 1 & 0 & 0 & 1 & 0 \\
0 & $-2$ \ & 1 & 1 & 0 & 0 & 0 & 1 \\
1 & 1 & $-2$ \ & $-2$ & 0 & 1 & 0 & 0 \\
1 & 1 & $-2$ & $-2$ \ & 0 & 1 & 0 & 0 \\
0 & 0 & 0 & 0 & $-2$ \ & 1 & 0 & 0 \\
0 & 0 & 1 & 1 & 1 & $-2$ \ & 0 & 0\\
1 & 0 & 0 & 0 & 0 & 0 & $-2$ \ & 0 \\
0 & 1 & 0 & 0 & 0 & 0 & 0 & $-2$ \
\end{array}
\right| \left|
\begin{array}{cccccccc}
$-2$ \ & 0 & 0 & 0 & 0 & 0 & 0 & 1 \\
0 & $-2$ \ & 0 & 1 & 0 & 1 & 1 & 0 \\
0 & 0 & $-2$ \ & 1 & 1 & 0 & 1 & 0 \\
0 & 1 & 1 & $-2$ \ & 0 & 0 & $-2$ & 1 \\
0 & 0 & 1 & 0 & $-2$ \ & 0 & 0 & 0 \\
0 & 1 & 0 & 0 & 0 & $-2$ \ & 0 & 0\\
0 & 1 & 1 & $-2$ & 0 & 0 & $-2$ \ & 1 \\
1 & 0 & 0 & 1 & 0 & 0 & 1 & $-2$ \
\end{array}
\right|\label{matp8-2}
\end{equation}

Our reduced program witnesses that the pairs of values of  invariants \ $(\mbox{\rm Card}, \mbox{\rm Ind}) $ \ for the fifth, sixth and seventh elements of the standard scale of the first of these morsifications (i.e., for virtual morsifications having respectively 4, 5 and 6 negative critical values) are equal to $(1216,-2), $ $(156,-3), $ and $(1216,-2)$. Similarly, the third, fourth, and fifth elements of the standard scale of the second virtual morsification (\ref{matp8-2}) (having respectively two, three and four negative critical values) have pairs of invariants equal to $(1216,-2),$ $ (258,-3), $ and $(1216,-2)$.  So, we get at least three new virtual components: one with \ $\mbox{\rm Ind} = -2$ \ and \ $\mbox{\rm Card} = 1216$, and two with $\mbox{\rm Ind} = -3$ and \ $\mbox{\rm Card}$ \ equal to 156 and 258. They fill in the gap, and statement 2 of Proposition \ref{proP2} is completely proved. Statement 3 follows immediately from Euler characteristic arguments, and statement 4 was proved during the consideration of  standard scales of morsifications (\ref{matp8}) and (\ref{matp8-3}).
 \hfill $\Box$

\end{document}